\documentclass[AMA,STIX1COL]{./ama/WileyNJD-v2}
\usepackage{moreverb}

\newcommand\BibTeX{{\rmfamily B\kern-.05em \textsc{i\kern-.025em b}\kern-.08em
T\kern-.1667em\lower.7ex\hbox{E}\kern-.125emX}}

\articletype{RESEARCH ARTICLE}%

\received{<day> <Month>, <year>}
\revised{<day> <Month>, <year>}
\accepted{<day> <Month>, <year>}

\usepackage{amsmath,amsfonts,amssymb}
\usepackage{mathtools}
\usepackage{siunitx}
\usepackage{tikz-dimline}

\usepackage{graphicx}
\usepackage{tikzexternal}

\usepackage{pgfplots, pgfplotstable}
\pgfplotstableset{col sep=comma, header=false}
\pgfplotsset{compat=1.16}
\pgfplotsset{
    every axis/.append style={
        thick,
        tick label style={font=\footnotesize},
        label style={font=\footnotesize},
        legend style={font=\footnotesize, line width=0.8pt},
    }
}

\definecolor{colorDiv0}{HTML}{762A83}
\definecolor{colorDiv1}{HTML}{9970AB}
\definecolor{colorDiv2}{HTML}{C2A5CF}
\definecolor{colorDiv3}{HTML}{A6DBA0}
\definecolor{colorDiv4}{HTML}{5AAE61}
\definecolor{colorDiv5}{HTML}{1B7837}

\definecolor{colorQ0}{HTML}{1F78B4}
\definecolor{colorQ1}{HTML}{FB9A99}
\definecolor{colorQ2}{HTML}{33A02C}
\definecolor{colorQ3}{HTML}{A6CEE3}
\definecolor{colorQ4}{HTML}{B2DF8A}

\pgfplotscreateplotcyclelist{myDivPlotCycleList}{%
    {colorDiv0}, 
    {colorDiv1}, 
    {colorDiv2}, 
    {colorDiv3}, 
    {colorDiv4}, 
    {colorDiv5}, 
}

\pgfplotscreateplotcyclelist{myQPlotCycleList}{%
    {colorQ1}, 
    {colorQ2}, 
    {colorQ3}, 
    {colorQ0}, 
    {colorQ4}, 
}

\usepackage{graphicx}
\usepackage[font=footnotesize,labelfont=bf,position=top, skip=0.5ex]{subcaption}

\usepackage{booktabs}
\usepackage{multirow}

\usepackage{rotating}

\def\*#1{\mathbf{#1}}
\def\nsd{{n_\mathrm{sd}}}
\def\nts{{n_\mathrm{ts}}}
\def\nex{{n_\mathrm{ex}}}

\def\ndf{{n_\mathrm{dof}}}

\def\Gt{\Gamma_t}

\newcommand{\ba}{\mathbf{a}}

\newcommand{\bx}{\mathbf{x}}

\newcommand{\bG}{\mathbf{G}}

\newcommand{\bI}{\mathbf{I}}
\newcommand{\bJ}{\mathbf{J}}

\newcommand{\bL}{\mathbf{L}}
\newcommand{\bM}{\mathbf{M}}

\newcommand{\zero}{\mathbf{0}}
\newcommand{\bxi}{{\boldsymbol{\xi}}}

\renewcommand{\div}{{\boldsymbol{\nabla}}}

\newcommand{\ddx}[2]{\frac{\partial #1}{\partial #2}}
\newcommand{\ddt}[3]{\left(\ddx{#1}{#2}\right)^{#3}}
\newcommand{\tri}[1]{$\bigtriangleup_{#1}$}

\newcommand{\prp}{piston ring pack }

\begin{document}

\title{Time-Continuous and Time-Discontinuous Space-Time Finite Elements for Advection-Diffusion Problems}

\author[1,2]{Max von Danwitz*}

\author[3]{Igor Voulis}

\author[2]{Norbert Hosters}

\author[2]{Marek Behr}

\authormark{M. v. DANWITZ \textsc{et al}}

\address[1]{\orgdiv{Institute for Mathematics and Computer-Based Simulation (IMCS)}, \orgname{University of the Bundeswehr Munich}, \orgaddress{\country{Germany}}}

\address[2]{\orgdiv{Chair for Computational Analysis of Technical Systems (CATS)}, \orgname{RWTH Aachen University}, \orgaddress{\country{Germany}}}

\address[3]{\orgdiv{Institute of Mathematics}, \orgname{Johannes Gutenberg University Mainz}, \orgaddress{\country{Germany}}}

\corres{*Max von Danwitz, Institute for Mathematics and Computer-Based Simulation (IMCS), University of the Bundeswehr Munich, Werner-Heisenberg-Weg 39, D-85577 Neubiberg, Germany. \email{max.danwitz@unibw.de}}

%\presentaddress{Present address}

\abstract[Abstract]{We construct four variants of space-time finite element discretizations based on linear tensor-product and simplex-type finite elements. The resulting discretizations are continuous in space, and continuous or discontinuous in time. In a first test run, all four methods are applied to a linear scalar advection-diffusion model problem. Then, the convergence properties of the time-discontinuous space-time finite element discretizations are studied in numerical experiments. Advection velocity and diffusion coefficient are varied, such that the parabolic case of pure diffusion (heat equation), as well as, the hyperbolic case of pure advection (transport equation) are included in the study. For each model parameter set, the $L_2$ error at the final time is computed for spatial and temporal element lengths ranging over several orders of magnitude to allow for an individual evaluation of the methods' spatial, temporal, and space-time accuracy. In the parabolic case, particular attention is paid to the influence of time-dependent boundary conditions. Key findings include a spatial accuracy of second order and a temporal accuracy between second and third order. The temporal accuracy tends towards third order depending on how advection-dominated the test case is, on the choice of the specific discretization method, and on the time-(in)dependence and treatment of the boundary conditions. Additionally, the potential of time-continuous simplex space-time finite elements for heat flux computations is demonstrated with a \prp test case.}


\jnlcitation{\cname{%
\author{v. Danwitz M}, 
\author{Voulis I}, 
\author{Hosters N}, and 
\author{Behr M}} (\cyear{2022}), 
\ctitle{Time-Continuous and Time-Discontinuous Space-Time Finite Elements for Advection-Diffusion Problems}, \cjournal{Journal}, \cvol{2022;00:x--x}.}

\maketitle

\footnotetext{Preprint submitted for publication}

%\title{Time-Continuous and Time-Discontinuous Space-Time Finite Elements for Advection-Diffusion Problems}
%%{\large \bf Time-Continuous and Time-Discontinuous Space-Time Finite Elements for Advection-Diffusion Problems} \par \vspace{6pt}
%\author{Max von Danwitz\footnote{Institute for Mathematics and Computer-Based Simulation (IMCS),
%Universit\"{a}t der Bundeswehr M\"{u}nchen, Werner-Heisenberg-Weg 39, 85577 Neubiberg,
%Germany. Email: max.danwitz@unibw.de} $^,$\footnotemark[3], Igor Voulis\footnote{Institute of Mathematics, Johannes Gutenberg University Mainz, Staudingerweg 9, 55128 Mainz, Germany. Email: ivoulis@uni-mainz.de}, Norbert Hosters\footnote{Chair for Computational
%Analysis of Technical Systems (CATS),
%Center for Simulation and Data Science (JARA-CSD),
%RWTH Aachen University, Schinkelstr. 2, 52062 Aachen, Germany. Email: \{danwitz, hosters, behr\}@cats.rwth-aachen.de}, Marek Behr\footnotemark[3]}
%\date{\today}
%%{\large Max von Danwitz\footnote{Chair for Computational
%%Analysis of Technical Systems (CATS),
%%RWTH Aachen University, 52056 Aachen,
%%Germany. Email: \{danwitz, behr\}@cats.rwth-aachen.de}, Igor Voulis\footnote{Institute of Mathematics, Johannes Gutenberg University Mainz, Staudingerweg 9, 55128 Mainz, Germany. Email: ivoulis@uni-mainz.de}, Marek Behr\footnotemark[1]} \par \vspace{6pt}

\maketitle

\begin{abstract}
\noindent

\end{abstract}

\section{Introduction}
\label{sec:intro}
\subsection{Motivation}

Multiple features make space-time finite elements an attractive solution strategy for time-dependent partial differential equations (PDE). First, space-time finite elements provide a uniform framework for error analysis as no distinction is made between spatial and temporal coordinates~\cite{hughes1988space}, which can also be used in adaptive refinement of the combined space-time mesh~\cite{langer2019coeff}. Moreover, space-time finite elements allow for parallel-in-time (PinT) computations which have inherently more potential for parallelization than spatial finite elements combined with a sequential time-stepping scheme~\cite{sivas2021air}. Furthermore, space-time finite elements are a natural choice to discretize time-dependent spatial computational domains, e.g., in fluid-structure interaction (FSI) simulations~\cite{tezduyar1992new, hubner2004monolithic, sathe2008modeling, spenke2020multi}. In particular, simplex space-time finite elements~\cite{behr2008simplex} can provide a boundary conforming space-time mesh for spatial domains that change topology over time~\cite{danwitz2021four}.

To benefit from these advantages, space-time finite elements have been used to perform simulations in various fields of computational fluid dynamics (CFD). Recent examples of simplex space-time simulations include the computation of complex fluid flows in production engineering applications~\cite{karyofylli2018simplex, karyofylli2019simplex} and the computation of dense granular flows~\cite{gesenhues2021simulating}. Likewise, compressible flows have been successfully simulated on unstructured space-time meshes~\cite{rendall2012conservative, wang2015high, danwitz2019simplex}. Note that the solution of transient three-dimensional problems with space-time finite elements requires four-dimensional meshes. Recent advances in generation~\cite{danwitz2021four, karabelas2019generating, boissonnat2021triangulating}, adaptation~\cite{caplan2020four}, and numerical handling~\cite{takizawa2019node,frontin2021foundations} of four-dimensional meshes mark the state-of-the-art in this active research field.

For efficiency considerations and refinement strategies, it is important to know the convergence behavior of the space-time finite element solution towards the physical or analytical solution of the simulated test case. However, for simulations based on the incompressible or compressible Naiver--Stokes equations it is an intricate task to estimate exact convergence orders, since numerical reference solutions can be influenced by round-off errors or implementation issues. Instead, we consider in this paper advection-diffusion problems \--- which lend themselves to an analytical solution \--- as a prototype for more complex flow problems~\cite{elman2014finite}. Based on the advection-diffusion equation, one can investigate the performance of numerical schemes with respect to transient, advective, and diffusive effects as well as their interplay. Besides, advection-diffusion equations also model a variety of physical problems, e.g., the concentration of a chemical species transported by an ambient flow or the temperature of a fluid streaming along a heated wall~\cite{elman2014finite}. Therefore, it is of great interest to analyze the convergence behavior of numerical schemes for advection-diffusion problems.

\subsection{Literature Review}

Shakib and Hughes~\cite{shakib1991new} present a Fourier analysis of space-finite elements with tensor-product structure applied to an advective-diffusive model problem with periodic boundary conditions. The method is found to be third order accurate with respect to the time step size for the pure advection and pure diffusion case. A summary of space-time finite element methods for convective transport problems is provided by Donea and Huerta along with numerical tests~\cite{donea2003finite}.

Moreover, linear tensor-product space-time finite elements can be related to a spatial discretization with finite elements and a temporal discretization with the Crank-Nicolson scheme~\cite{behr2008simplex}. Studies of this resulting method often focus either on parabolic problems (heat equation)~\cite{lozinski2009anisotropic} or on the pure advection case (transport equation)~\cite{dubuis2018adaptive}. Moreover, a Crank–Nicolson type space-time finite element method for evolution problems on moving meshes is proposed and analyzed by Hansbo~\cite{hansbo2000crank}. The method uses tensor product elements that are inclined in space-time with a slope given by the convection velocity. It is reported that the aligned space-time orientation improves the precision and facilitates the solution of the discrete system.

Focusing on the parabolic limit case, time-continuous tensor-product space-time finite elements have been analyzed by Aziz and Monk~\cite{aziz1989continuous}. In more recent works, also unstructured space-time finite elements which do not require any tensor-product structure are addressed, e.g., by Steinbach~\cite{steinbach2015space}. Furthermore, Langer and Schafelner~\cite{langer2019coeff, langer2019spacetime} investigate the scaling behavior of unstructured space-time finite element methods for parabolic problems in parallel computations. Note that this work is also extended to hexahedral space-time discretizations~\cite{langer2021hex}. Moreover, Langer and Zank propose and investigate new efficient direct solvers for time-continuous tensor-product discretizations of the parabolic initial boundary value problem~\cite{langer2020efficient}. The influence of linear constraints, e.g., time-dependent Dirichlet boundary conditions, on discontinuous Galerkin time discretization methods for parabolic problems is treated by Voulis and Reusken~\cite{voulis2019discontinuous}.

\subsection{Scientific Novelty and Limitations}
To the best of the authors' knowledge, there is no previous comprehensive numerical study that analyses the convergence behavior of tensor-product and simplex-type finite elements for the complete range of model parameters of advection-diffusion problems and for spatial and temporal element sizes over several orders of magnitude. On the one hand, the computational evaluation of the convergence behavior is advantageous in the sense that a simple variation of the input parameters allows to switch from a parabolic to a hyperbolic problem. Therefore, the computational approach facilitates a study of the precise influence of parameter variations. On the other hand, the numerical study is limited to specific test cases and for those considers only the $L_2$ norm and a nodal measure of the error at the final time. Still, it is expected that the results also hold for other test cases of similar nature. 

\subsection{Paper Organization}

In the remainder, we proceed as follows. In Section~\ref{sec:method}, four space-time discretizations are presented and descriptive naming is proposed. In Section~\ref{sec:adEquation}, we apply the methods to an initial boundary value problem based on the advection-diffusion equation. Section~\ref{sec:convergence} collects the results of a computational error analysis of the time-discontinuous discretizations and compares the results with the theoretically expected convergence behavior. In Section~\ref{sec:prp}, we demonstrate the particular potential of simulations on time-continuous simplex space-time meshes in a \prp application. Concluding remarks are offered in Section~\ref{sec:conclusion}.

\section{Method Classification}
\label{sec:method}
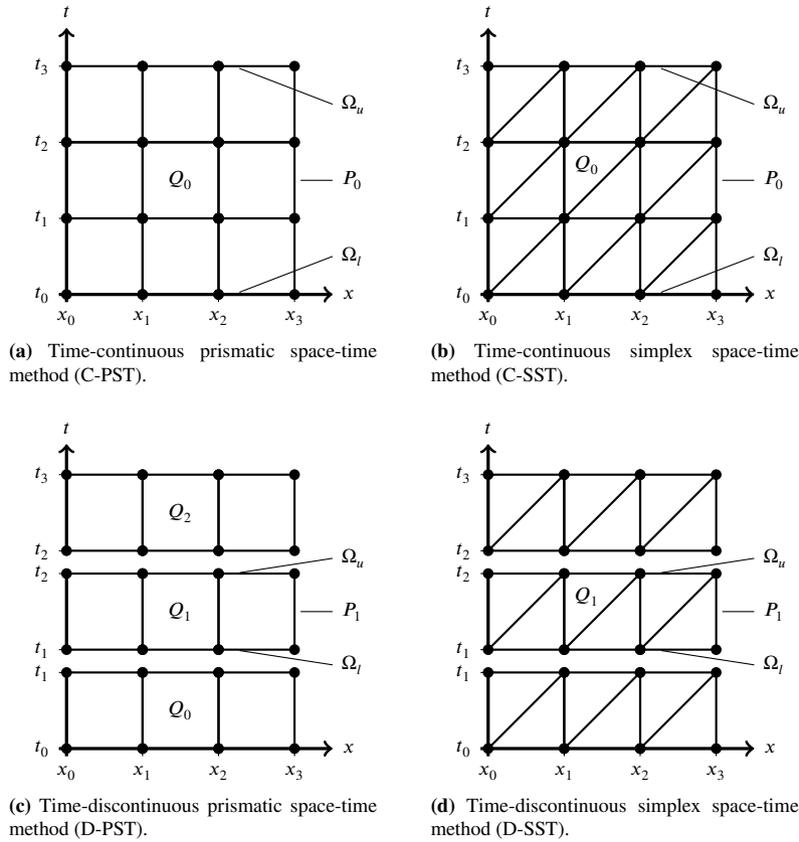
\begin{figure}
\captionsetup[sub]{position=bottom}
\centering
\subcaptionbox{Time-continuous prismatic space-time method (C-PST).\label{fig:CPST}}{\newcommand*\cols{2}
\newcommand*\rows{2}
%\tikzsetnextfilename{spacetime-cpst}
\begin{tikzpicture}[
    scale=1,
    axis/.style={very thick, ->},
    important line/.style={thick},
    dashed line/.style={dashed, thin},
    pile/.style={thick, ->, shorten <=2pt, shorten
    >=2pt},
    every node/.style={color=black},
    pointer/.style={shorten <=0.08cm},
    ]
    %x-axis
    \draw[axis] (-0.1,0)  -- (3.5,0) node(xline)[right]{$x$};
    \foreach \x in {0,1,2,3}
    \draw (\x cm,0.1) -- (\x cm,-0.1) node[anchor=north] {$x_{\x}$};
    % t-axis
    \draw[axis] (0,-0.1) -- (0,3.5) node(yline)[above]{$t$};
    \draw (0.1,0cm) -- (-0.1, 0cm) node[anchor=east] {$t_{0}$};
    \draw (0.1,1cm) -- (-0.1, 1cm) node[anchor=east] {$t_{1}$};
    \draw (0.1,2.0cm) -- (-0.1, 2.0cm) node[anchor=east] {$t_{2}$};
    \draw (0.1,3.0cm) -- (-0.1, 3.0cm) node[anchor=east] {$t_{3}$};
    % Lines
    \foreach \row in {0, 1, ...,\rows} {
     		\foreach \col in {0, 1, ...,\cols} {
   				 \draw[important line] (\col,\row) coordinate (A)  -- (\col+1,\row) coordinate (B);
   				 \draw[important line] (\col,\row+1) coordinate (A)  -- (\col+1,\row+1) coordinate (B);
   				 \draw[important line] (\col,\row) coordinate (A)  -- (\col,\row+1) coordinate (B);
   				 \fill (A) circle (2pt);
   				 \fill (B) circle (2pt);
   				  \draw[important line] (\col+1,\row) coordinate (A)  -- (\col+1,\row+1) coordinate (B);
   				 \fill (A) circle (2pt);
   				 \fill (B) circle (2pt);
   		 }
   	}
   \node (Q) at (1.5,1.5) {$Q_0$};
   \draw[pointer] (3,1.5cm) -- (3.5, 1.5) node[anchor=west] {$P_{0}$};
   \draw[pointer] (2.2,0cm) -- (3.5, 0.5) node[anchor=west] {$\Omega_l$};
        \draw[pointer] (2.2,3cm) -- (3.5, 2.5) node[anchor=west] {$\Omega_u$};
\end{tikzpicture}}\qquad
\subcaptionbox{Time-continuous simplex space-time method (C-SST).\label{fig:CSST}}{\newcommand*\cols{2}
\newcommand*\rows{2}
%\tikzsetnextfilename{spacetime-csst}
\begin{tikzpicture}[
    scale=1,
    axis/.style={very thick, ->},
    important line/.style={thick},
    pointer/.style={shorten <=0.08cm},
    every node/.style={color=black}
    ]
    %x-axis
    \draw[axis] (-0.1,0)  -- (3.5,0) node(xline)[right]{$x$};
    \foreach \x in {0,1,2,3}
    \draw (\x cm,0.1) -- (\x cm,-0.1) node[anchor=north] {$x_{\x}$};
    % t-axis
%    \draw[very thick] (0,-0.1) -- (0,1);
%    \draw[very thick] (0,1.3) -- (0,2.3);
    \draw[axis] (0,-0.1) -- (0,3.5) node(yline)[above]{$t$};
    \draw (0.1,0cm) -- (-0.1, 0cm) node[anchor=east] {$t_{0}$};
    \draw (0.1,1cm) -- (-0.1, 1cm) node[anchor=east] {$t_{1}$};
    \draw (0.1,2.0cm) -- (-0.1, 2.0cm) node[anchor=east] {$t_{2}$};
    \draw (0.1,3.0cm) -- (-0.1, 3.0cm) node[anchor=east] {$t_{3}$};
    % Lines
    \foreach \row in {0, 1, ...,\rows} {
     		\foreach \col in {0, 1, ...,\cols} {
   				 \draw[important line] (\col,\row) coordinate (A)  -- (\col+1,\row) coordinate (B);
   				 \draw[important line] (\col,\row+1) coordinate (A)  -- (\col+1,\row+1) coordinate (B);
   				 \draw[important line] (\col,\row) coordinate (A)  -- (\col,\row+1) coordinate (B);
   				 \fill (A) circle (2pt);
   				 \fill (B) circle (2pt);
   				  \draw[important line] (\col+1,\row) coordinate (A)  -- (\col+1,\row+1) coordinate (B);
   				 \fill (A) circle (2pt);
   				 \fill (B) circle (2pt);
				   \draw[important line] (\col,\row) coordinate (A)  -- (\col+1,\row+1) coordinate (B);
   		 }
   	}
   \node (Q) at (1.3,1.7) {$Q_0$};
   \draw[pointer] (3,1.5cm) -- (3.5, 1.5) node[anchor=west] {$P_{0}$};
     \draw[pointer] (2.2,0cm) -- (3.5, 0.5) node[anchor=west] {$\Omega_l$};
        \draw[pointer] (2.2,3cm) -- (3.5, 2.5) node[anchor=west] {$\Omega_u$};
\end{tikzpicture}}\\ \vspace{2ex}
\subcaptionbox{Time-discontinuous prismatic space-time method (D-PST).\label{fig:DPST}}{\newcommand*\cols{2}
\newcommand*\rows{3}

%\tikzsetnextfilename{spacetime-dpst}
\begin{tikzpicture}[
    scale=1,
    axis/.style={very thick, ->},
    important line/.style={thick},
    pointer/.style={shorten <=0.08cm},
    every node/.style={color=black}
    ]
    %x-axis
    \draw[axis] (-0.1,0)  -- (3.5,0) node(xline)[right]{$x$};
    \foreach \x in {0,1,2,3}
    \draw (\x cm,0.1) -- (\x cm,-0.1) node[anchor=north] {$x_{\x}$};
    % t-axis
    \draw[very thick] (0,-0.1) -- (0,1);
    \draw[very thick] (0,1.3) -- (0,2.3);
    \draw[axis] (0,2.6) -- (0,4) node(yline)[above]{$t$};
    \draw (0.1,0cm) -- (-0.1, 0cm) node[anchor=east] {$t_{0}$};
    \draw (0.1,1cm) -- (-0.1, 1cm) node[anchor=east] {$t_{1}$};
    \draw (0.1,1.3cm) -- (-0.1, 1.3cm) node[anchor=east] {$t_{1}$};
    \draw (0.1,2.3cm) -- (-0.1, 2.3cm) node[anchor=east] {$t_{2}$};
    \draw (0.1,2.6cm) -- (-0.1, 2.6cm) node[anchor=east] {$t_{2}$};
    \draw (0.1,3.6cm) -- (-0.1, 3.6cm) node[anchor=east] {$t_{3}$};
    % Lines
    \foreach \row in {0, 1.3, ...,\rows} {
     		\foreach \col in {0, 1, ...,\cols} {
   				 \draw[important line] (\col,\row) coordinate (A)  -- (\col+1,\row) coordinate (B);
   				 \draw[important line] (\col,\row+1) coordinate (A)  -- (\col+1,\row+1) coordinate (B);
   				 \draw[important line] (\col,\row) coordinate (A)  -- (\col,\row+1) coordinate (B);
   				 \fill (A) circle (2pt);
   				 \fill (B) circle (2pt);
   				  \draw[important line] (\col+1,\row) coordinate (A)  -- (\col+1,\row+1) coordinate (B);
   				 \fill (A) circle (2pt);
   				 \fill (B) circle (2pt);
   		 }
   	}
	\node (Q) at (1.5,0.5) {$Q_0$};
	\node (Q) at (1.5,1.8) {$Q_1$};
	\node (Q) at (1.5,3.1) {$Q_2$};
   	\draw[pointer] (3,1.8cm) -- (3.5, 1.8) node[anchor=west] {$P_{1}$};
     	\draw[pointer] (2.2,1.3cm) -- (3.5, 1.1) node[anchor=west] {$\Omega_l$};
        \draw[pointer] (2.2,2.3cm) -- (3.5, 2.5) node[anchor=west] {$\Omega_u$};
\end{tikzpicture}}\qquad
\subcaptionbox{Time-discontinuous simplex space-time method (D-SST).\label{fig:DSST}}{\newcommand*\cols{2}
\newcommand*\rows{3}
%\tikzsetnextfilename{spacetime-dsst}
\begin{tikzpicture}[
    scale=1,
    axis/.style={very thick, ->},
    important line/.style={thick},
    newEdge/.style={thick, ->, >=circle},
        pointer/.style={shorten <=0.08cm},
    every node/.style={color=black}
    ]
    %x-axis
    \draw[axis] (-0.1,0)  -- (3.5,0) node(xline)[right]{$x$};
    \foreach \x in {0,1,2,3}
    \draw (\x cm,0.1) -- (\x cm,-0.1) node[anchor=north] {$x_{\x}$};
    % t-axis
    \draw[very thick] (0,-0.1) -- (0,1);
       \draw[very thick] (0,1.3) -- (0,2.3);
    \draw[axis] (0,2.6) -- (0,4) node(yline)[above]{$t$};
    \draw (0.1,0cm) -- (-0.1, 0cm) node[anchor=east] {$t_{0}$};
    \draw (0.1,1cm) -- (-0.1, 1cm) node[anchor=east] {$t_{1}$};
    \draw (0.1,1.3cm) -- (-0.1, 1.3cm) node[anchor=east] {$t_{1}$};
    \draw (0.1,2.3cm) -- (-0.1, 2.3cm) node[anchor=east] {$t_{2}$};
        \draw (0.1,2.6cm) -- (-0.1, 2.6cm) node[anchor=east] {$t_{2}$};
            \draw (0.1,3.6cm) -- (-0.1, 3.6cm) node[anchor=east] {$t_{3}$};
    % Lines
    \foreach \row in {0, 1.3, ...,\rows} {
     		\foreach \col in {0, 1, ...,\cols} {
   				 \draw[important line] (\col,\row) coordinate (A)  -- (\col+1,\row) coordinate (B);
   				 \draw[important line] (\col,\row+1) coordinate (A)  -- (\col+1,\row+1) coordinate (B);
				 \draw[important line] (\col,\row) coordinate (A)  -- (\col,\row+1) coordinate (B);
   				 \fill (A) circle (2pt);
   				 \fill (B) circle (2pt);
   				  \draw[important line] (\col+1,\row) coordinate (A)  -- (\col+1,\row+1) coordinate (B);
   				 \fill (A) circle (2pt);
   				 \fill (B) circle (2pt);
				  \draw[important line] (\col,\row) coordinate (A)  -- (\col+1,\row+1) coordinate (B);
   		 }
   		    	%\draw[important line] (0,\row) coordinate (A)  -- (1,\row+1) coordinate (B);
%			\draw[important line] (0,\row+1) coordinate (A)  -- (1,\row) coordinate (B);
%   				 \fill (A) circle (2pt);
%   				 \fill (B) circle (2pt);
%   				 
%   				  \draw[important line] (1,\row) coordinate (A)  -- (2,\row+0.5) coordinate (B);
%   				   \fill (A) circle (2pt);
%   				  \draw[important line] (1,\row+1) coordinate (A)  -- (2,\row+0.5) coordinate (B);
%   				  \fill (A) circle (2pt);
%   				  \draw[important line] (2,\row+0.5) coordinate (A)  -- (3,\row) coordinate (B);
%   				   \fill (A) circle (2pt);
%   				   \fill (B) circle (2pt);
%   				  \draw[important line] (2,\row+0.5) coordinate (A)  -- (3,\row+1) coordinate (B);
%   				  \fill (A) circle (2pt);   	
%   				  \fill (B) circle (2pt);			  
   				  
   	}
   	   \node (Q) at (1.3,2.0) {$Q_1$};
   \draw[pointer] (3,1.8cm) -- (3.5, 1.8) node[anchor=west] {$P_{1}$};
     \draw[pointer] (2.2,1.3cm) -- (3.5, 1.1) node[anchor=west] {$\Omega_l$};
        \draw[pointer] (2.2,2.3cm) -- (3.5, 2.5) node[anchor=west] {$\Omega_u$};
\end{tikzpicture}}\\
\caption{Space-time discretization methods.}
\label{fig:ST}
\end{figure}

To introduce the specific space-time discretizations investigated in this work, the naming of involved entities is briefly reviewed~\cite{hughes1988space, danwitz2019simplex}. We consider a spatial computational domain $\Omega \subset \mathbb{R}^{\nsd}$, where $\nsd$ denotes the number of spatial dimensions. That domain $\Omega$ and a time interval, $I = [0,t_f] \subset \mathbb{R}$, span the space-time continuum $Q \subset \mathbb{R}^{\nsd+1}$. In the following, we consider four ways to approximate the solution of PDEs on $Q$ with finite elements. Sample slicings $Q^h$ of the space-time domain $Q =[x_0, x_3]\times [t_0,t_3]$ are shown in Figure~\ref{fig:ST}. For the sake of clarity, the spatial domain $\Omega$ remains constant over time in these drawings. However, the proposed methods can also be applied to time-dependent spatial domains $\Omega(t)$~\cite{tezduyar1992new, danwitz2021four, karabelas2019generating}. The first two discretization techniques (Figure~\ref{fig:CPST} and~\ref{fig:CSST}) seek an approximation that is continuous across $Q$. In contrast, the second two (Figure~\ref{fig:DPST} and~\ref{fig:DSST}) seek an approximation that is discontinuous at certain times, which leads to a discontinuous Galerkin method for the temporal discretization. In these time-discontinuous cases, $Q$ is sliced into space-time slabs $Q_n$. As indicated in the drawings of Figure~\ref{fig:DPST} and~\ref{fig:DSST}, the boundary of each space-time discretization consists of three parts: the spatial discretization at the lower time level $\Omega^h_{\text{l}} = \Omega^h(t=t_n)$, the spatial discretization at the upper time level $\Omega^h_{\text{u}}= \Omega^h(t=t_{n+1})$, and the discretization of the space-time boundary $P \subset \mathbb{R}^{\nsd}$ which is the temporal evolution of the spatial domain boundary $\Gamma \subset \mathbb{R}^{\nsd-1}$. The size of space-time slabs in temporal direction is denoted by $\Delta t$. To later apply one uniform finite element formulation for the time-continuous and time-discontinuous cases, we regard the complete space-time domain $Q$ in the time-continuous case as space-time slab $Q_0$.

\begin{figure}
\captionsetup[sub]{position=bottom}
\centering
\subcaptionbox{Local temporal refinement in D-SST mesh.\label{fig:SST}}{\newcommand*\cols{2}
\newcommand*\rows{3}
%\tikzsetnextfilename{spacetime-sst}
\begin{tikzpicture}[
    scale=1,
    axis/.style={very thick, ->, >=stealth'},
    important line/.style={thick},
    newEdge/.style={thick, ->, >=circle},
            pointer/.style={shorten <=0.08cm},
    every node/.style={color=black}
    ]
    %x-axis
    \draw[axis] (-0.1,0)  -- (3.5,0) node(xline)[right]{$x$};
    \foreach \x in {0,1,2,3}
    \draw (\x cm,0.1) -- (\x cm,-0.1) node[anchor=north] {$x_{\x}$};
    % t-axis
    \draw[very thick] (0,-0.1) -- (0,1);
       \draw[very thick] (0,1.3) -- (0,2.3);
    \draw[axis] (0,2.6) -- (0,4) node(yline)[above]{$t$};
    \draw (0.1,0cm) -- (-0.1, 0cm) node[anchor=east] {$t_{0}$};
    \draw (0.1,1cm) -- (-0.1, 1cm) node[anchor=east] {$t_{1}$};
    \draw (0.1,1.3cm) -- (-0.1, 1.3cm) node[anchor=east] {$t_{1}$};
    \draw (0.1,2.3cm) -- (-0.1, 2.3cm) node[anchor=east] {$t_{2}$};
        \draw (0.1,2.6cm) -- (-0.1, 2.6cm) node[anchor=east] {$t_{2}$};
            \draw (0.1,3.6cm) -- (-0.1, 3.6cm) node[anchor=east] {$t_{3}$};
    % Lines
    \foreach \row in {0, 1.3, ...,\rows} {
     		\foreach \col in {0, 1, ...,\cols} {
   				 \draw[important line] (\col,\row) coordinate (A)  -- (\col+1,\row) coordinate (B);
   				 \draw[important line] (\col,\row+1) coordinate (A)  -- (\col+1,\row+1) coordinate (B);
   				 \draw[important line] (\col,\row) coordinate (A)  -- (\col,\row+1) coordinate (B);
   				 \fill[black] (A) circle (2pt);
   				 \fill[black] (B) circle (2pt);
   				  \draw[important line] (\col+1,\row) coordinate (A)  -- (\col+1,\row+1) coordinate (B);
   				 \fill[black] (A) circle (2pt);
   				 \fill[black] (B) circle (2pt);
   		 }
   		    	\draw[important line] (0,\row) coordinate (A)  -- (1,\row+1) coordinate (B);
   				 \fill[black] (A) circle (2pt);
   				 \fill[black] (B) circle (2pt);
   				 
   				  \draw[important line] (1,\row) coordinate (A)  -- (2,\row+0.5) coordinate (B);
   				   \fill[black] (A) circle (2pt);
   				  \draw[important line] (1,\row+1) coordinate (A)  -- (2,\row+0.5) coordinate (B);
   				  \fill[black] (A) circle (2pt);
   				  \draw[important line] (2,\row+0.5) coordinate (A)  -- (3,\row) coordinate (B);
   				   \fill[black] (A) circle (2pt);
   				   \fill[black] (B) circle (2pt);
   				  \draw[important line] (2,\row+0.5) coordinate (A)  -- (3,\row+1) coordinate (B);
   				  \fill[black] (A) circle (2pt);   	
   				  \fill[black] (B) circle (2pt);			  
   				  
   	}
   	   \node (Q) at (1.5,1.8) {$Q_1$};
   \draw[pointer] (3,1.8cm) -- (3.5, 1.8) node[anchor=west] {$P_{1}$};
     \draw[pointer] (2.2,1.3cm) -- (3.5, 1.1) node[anchor=west] {$\Omega_l$};
        \draw[pointer] (2.2,2.3cm) -- (3.5, 2.5) node[anchor=west] {$\Omega_u$};
\end{tikzpicture}}\qquad
\subcaptionbox{Fully-unstructured C-SST mesh.\label{fig:UST}}{%\tikzsetnextfilename{spacetime-ust}
\begin{tikzpicture}[
    axis/.style={very thick, ->, >=stealth'},
    important line/.style={thick},
    pile/.style={thick, ->, >=stealth', shorten <=2pt, shorten
    >=2pt},
    pointer/.style={shorten <=0.08cm},
    every node/.style={color=black}
    ]
         %x-axis
    \draw[axis] (-0.1,0)  -- (3.5,0) node(xline)[right]{$x$};
    \foreach \x in {0,1,2,3}
    \draw (\x ,0.1) -- (\x,-0.1) node[anchor=north] {$x_{\x}$};
    % t-axis
    \draw[axis] (0,-0.1) -- (0,4) node(yline)[above]{$t$};
    \draw (0.1,0cm) -- (-0.1, 0cm) node[anchor=east] {$t_{0}$};
    \draw (0.1,1cm) -- (-0.1, 1cm) node[anchor=east] {$t_{1}$};
    \draw (0.1,2cm) -- (-0.1, 2cm) node[anchor=east] {$t_{2}$};
    \draw (0.1,3cm) -- (-0.1, 3cm) node[anchor=east] {$t_{3}$};
\begin{axis}[axis lines=none,
	xmin = 0,
	xmax = 3,
	ymin =0,
	ymax = 3,
	x=1cm,
	y=1cm
	]
	\addplot [patch, patch refines=0, mesh, thick, black, patch table={spacetime/d1_triangles.dat}] table {spacetime/d1_points.dat}; 
  	\addplot [mark=*, draw=black, fill=black, only marks] table {spacetime/d1_points.dat}; 
\end{axis}
  	\node (Q) at (1.5,1.5) {$Q_0$};
   	\draw[pointer] (3,1.5cm) -- (3.5, 1.5) node[anchor=west] {$P_{0}$};
     	\draw[pointer] (2.2,0cm) -- (3.5, 0.5) node[anchor=west] {$\Omega_l$};
        \draw[pointer] (2.2,3cm) -- (3.5, 2.5) node[anchor=west] {$\Omega_u$};
        % Extra lines
        \draw[thick] (3.0, 0.0) -- (3.0, 3.0);
        \draw[thick] (0.0, 3.0) -- (3.0, 3.0);
\end{tikzpicture}}
\caption{Variants of simplex space-time discretization methods.}
\label{fig:stAlternative}
\end{figure}
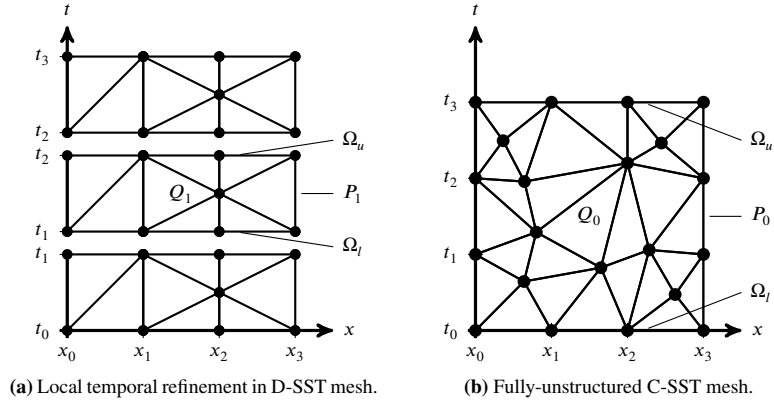

Both time-discretization approaches can be combined with prismatic elements with tensor-product structure, or simplex elements. The combinations form the four discretization methods C-PST, C-SST, D-PST, and D-SST. In PST methods, a discretization of $Q_n$ with prismatic space-time elements can be easily obtained by extrusion of a spatial discretization of $\Omega$ in time. C-PST is a continuous finite element discretization in space and time as described by Aziz and Monk~\cite{aziz1989continuous}. When combined with linear shape functions, it is also known as cg(1)cg(1). The time-discontinuous D-PST method is also referred to as cg(1)dg(1) for example by Quarteroni et al.~\cite{quarteroni2006numerical}. SST discretizations can be generated by subdividing the prismatic elements into simplex elements $Q^e_n$ (Figure~\ref{fig:CSST} and~\ref{fig:DSST}). More complex SST mesh generation procedures also allow for local temporal refinement by node insertion~\cite{behr2008simplex} or fully-unstructured space-time meshes~\cite{danwitz2019simplex} as shown in Figure~\ref{fig:stAlternative}.

For each space-time slab $Q_n$, an $H^1$-conformal finite element approximation space $H^1_{h,n}$ is constructed based on one of the discussed discretizations and element basis functions~\cite{danwitz2019simplex}. In case of PST discretizations, we consider the $\mathbb{PR}_1$ and $\mathbb{Q}_1$ basis functions of the simplex-based prismatic and cuboid Lagrange finite element. In case of SST discretizations, we use the $\mathbb{P}_1$ basis functions of the simplical Lagrange finite element ~\cite{ern2013theory}. The four space-time discretizations introduced, are now employed in the solution of advection-diffusion problems.

\section{Application to Advection-Diffusion Equation}
\label{sec:adEquation}

We consider the time-dependent linear advection-diffusion equation
\begin{equation}
\text{res}(u) \coloneqq  \ddx{u}{t} + \ba \cdot \nabla u - k \, \Delta u = 0. \label{eq:ADstrong}
\end{equation}
Therein, the scalar unknown, $u(\bx, t)$ is a function of the spatial coordinates ($\bx = \left(x,y,z\right)^T$ for $\nsd=3$) and time. The advection velocity is a given vector $\ba$, and the diffusion coefficient is denoted by $k$. As usual, the Laplacian of $u$ abbreviates $\Delta u = \nabla \cdot \nabla u$, based on the spatial gradient $\nabla u$. Advection velocity and diffusion coefficient can be varied, such that the parabolic case of pure diffusion ($\ba = \zero$), as well as the hyperbolic case of pure advection ($k=0$) are included. In the former case, Equation~\eqref{eq:ADstrong} is the heat equation, in the latter case the transport equation. Furthermore, the above equation lends itself to an analytic solution, hence, facilitating a computational error analysis as presented in Section~\ref{sec:convergence}.

A general characterization of advection-diffusion problems can be achieved with the dimensionless P\'{e}clet number 
\newcommand{\Peclet}{\mathrm{Pe}}
\begin{equation} \Peclet \coloneqq L_c \frac{a}{k}. \label{eq:Peclet} \end{equation}
Therein, a scalar measure of the advection speed, $a= \|\ba\|$, is related to the diffusion coefficient $k$ and scaled by a characteristic length $L_c$. As the dimensionless number compares the importance of advective and diffusive effects for a given test case, one can typically expect solutions with smaller gradients for test cases with lower P\'{e}clet number (when diffusion dominates).

To construct an initial boundary value problem, let us consider again a computational domain $Q$ as for example shown in Figure~\ref{fig:UST}. The associated space-time boundary $P$ is assumed to consist of a Dirichlet part $P^D$ and a Neumann part $P^N$, such that $P = P^D \cup P^N$  and $ P^D \cap P^N = \emptyset$. Then, we obtain an initial boundary value problem, as we require Equation~\eqref{eq:ADstrong} to hold on $Q$, along with a known initial condition $u_0$ and given Dirichlet boundary conditions $g$. The complete statement of the initial boundary value problem reads
\begin{align}
\label{eq:IBVP}
\text{IBVP} \quad
\begin{cases}
	\text{res}(u(\bx, t)) = 0, &\;  \text{on} \; Q,  \\
	u(\bx, t) = u_0(\bx),& \;  \text{at} \; t=0,   \\
	u(\bx, t) = g(\bx, t) , &\;  \text{on} \; P^D. \\
\end{cases}
\end{align}

When applying one of the discretization techniques described in Section~\ref{sec:method} to $Q$, the initial condition is enforced on $\Omega_l$ of the space-time slab $Q_0$. The part of a space-time slab $Q_n$, where Dirichlet boundary conditions are prescribed is denoted by $P_n^D$. A suitable interpolation of the Dirichlet boundary data $g^h$ allows us to define the trial function space
\begin{equation}
\mathcal{S}_{h,n} = \left\lbrace u^h \in  H^1_{h,n}\; \Bigg| \; u^h = g^h \; \text{on}\; P^D_n \right\rbrace 
\end{equation}
and the test function space
\begin{equation}
\mathcal{V}_{h,n} = \left\lbrace w^h \in H^1_{h,n} \; \Bigg| \; w^h = 0  \;\text{on}\; P^D_n \right\rbrace.
\end{equation}
Considering that for time-discontinuous discretization methods the finite element approximation is discontinuous at the space-time slab boundaries $\Omega_l$ and $\Omega_u$, let $\left(u^h\right)^{\pm}_n$ abbreviate $\lim_{\varepsilon \to 0} u^h(t_n \pm \varepsilon)$.

Using these definitions, a discretized weak form of the initial boundary value problem can be stated as follows:
For given initial conditions $\left(u^h\right)^-_0 = u^h_0$, find $u^h \in \mathcal{S}_{h,n}$ such that on each time slab $Q_n$ and for all $w^h \in \mathcal{V}_{h,n} $
%\begin{align}
\begin{alignat}{2}
\label{eq:weakAD}
0 =& & &\int_{Q_n} w^h \cdot \left( \ddx{u^h}{t}  +  \ba \cdot \nabla u^h \right) dQ \\ \nonumber
		    &+ &  &\int_{Q_n} \nabla w^h \cdot \left( k \nabla u^h \right) dQ \\ \nonumber
		   &+ & &\int_{\Omega_{l}} \left(w^h\right)^+_n \cdot \left[ \left(u^h \right)^+_n - \left(u^h\right)^-_n \right] d\Omega \\ \nonumber
		   &+ & &\int_{Q_n} \! \left( \ddx{w^h}{t} + \ba \cdot \nabla w^h \right) \cdot \tau_\mathrm{SUPG}  \cdot \text{res}\left(u^h\right)  dQ . \nonumber
%\end{align}
\end{alignat}
In the weak form above, the diffusion term was modified using integration by parts. The resulting boundary integral vanishes, since the test functions vanish on $P_n^D$ and homogeneous Neumann boundary conditions are assumed on $P_n^N$. Moreover, the initial condition as well as the continuity of $u^h$ between time slabs is weakly enforced with the integral over the spatial computational domain $\Omega_l$, the so-called jump term. The stability of the formulation is achieved with a SUPG term in the fourth integral~\cite{brooks1982streamline}. We define the stabilization parameter $\tau_\mathrm{SUPG}$ as 
\begin{equation}
\label{eq:tauAD}
\tau_{\mathrm{SUPG}} = \left( \left[  \begin{array}{c} \ba \\ 1 \\ \end{array} \right] \cdot \bG \left[  \begin{array}{c} \ba \\ 1 \\ \end{array} \right]+  \left( C_{\mathrm{inv}}\frac{ k }{h_{s}^2} \right)^2 \right)^{-\frac{1}{2}},
\end{equation}
which accounts for local characteristics of the initial boundary value problem. In the first term, the space-time element metric $\bG$ is used to include directional element length information. The metric tensor, 
\begin{equation}
\bG = \ddt{(\bxi,\tau)}{(\bx,t)}{T} \bM \,  \left(\ddx{(\bxi,\tau)}{(\bx,t)}\right), \label{eq:cov}
\end{equation}
is based on the inverse of the Jacobian associated with the mapping from reference coordinates,~$(\bxi,\tau)$, to physical coordinates, $(\bx,t)$. Moreover, the metric tensor includes a square matrix $\bM$ of size $\nsd+1$, which accounts for the mapping to a regular reference element counteracting the influence of the element's node numbering~\cite{danwitz2019simplex}. A further analysis of node-numbering invariant element length measures for simplex elements is presented by Takizawa et al.~\cite{takizawa2019node}. Explicit forms of $\bM$ for simplex elements read for $d=2,3,4$, respectively,
\begin{align}
\label{eq:Md}
	\bM_{d=2} = \frac{1}{\sqrt{3}} \left( \begin{array}{cc}
	2 & 1  \\
	1 & 2  \\
	\end{array} \right), \quad%
		\bM_{d=3} = \frac{1}{\sqrt[3]{4}} \left( \begin{array}{ccc}
	2 & 1 & 1  \\
	1 & 2 & 1 \\
	1 & 1 & 2 \\
	\end{array} \right), \quad
			\bM_{d=4} = \frac{1}{\sqrt[4]{5}} \left( \begin{array}{cccc}
	2 & 1 & 1  & 1\\
	1 & 2 & 1 & 1\\
	1 & 1 & 2 & 1\\
	1 & 1 & 1 &2 \\
	\end{array} \right).
\end{align}

For other element types, an appropriate matrix $\bM$ is substituted. Recalling that $\bM$ accounts for the mapping to a regular reference element, it is clear that discretizations with pure tensor-product reference elements ($\mathbb{Q}_1$) do not need an additional mapping\textemdash as the reference element is already regular. Therefore, $\bM$ can simply be replaced by the identity matrix. In case of a simplex-based prismatic reference element ($\mathbb{PR}_1)$, the partial tensor-product structure of the reference element is reflected in the choice of $\bM$ as shown below
 
\begin{equation}
\label{eq:Mspec}
\bM =
\begin{cases}
\quad \bM_{d=\nsd+1}  & \mathbb{P}_1\\
	\\
\quad \bI &  \mathbb{Q}_1\\
	\\
\left( \begin{array}{cc}
    	\bM_{d=\nsd} & \zero \\
	\zero^T  &1 \\
\end{array} \right) & \mathbb{PR}_1. \\
\end{cases}	
\end{equation}

In the second term of Equation~\eqref{eq:tauAD}, the diffusive contribution to $\tau_{\mathrm{SUPG}}$ requires a measure of the spatial element length $h_s$.  For all considered element types, the length $h_s$ is obtained from the spatial part of the metric tensor $\bG_s = \left[ \bG \right]_{\nsd \times \nsd}$ as
\begin{equation} \frac{1}{{h_s}^2} = \sqrt{\bG_s \colon \bG_s}, \end{equation}
where the colon operator denotes the double contraction $\bG \colon \bG = \sum_{i,j} G_{ij} \cdot G_{ij}$. Moreover, the constant $C_{\mathrm{inv}}$ scales the diffusive contribution to $\tau_{\mathrm{SUPG}}$. Inspired by an inverse estimate inequality proven in~\cite{knechtges2018simulation}, we chose for $ \mathbb{P}_1$ and $\mathbb{PR}_1$ discretizations
\begin{equation}
C_{\mathrm{inv}} = (\nsd+1)^2(\nsd+2) =
\begin{cases}
12 & \nsd = 1, \\
36 & \nsd = 2, \\
80 & \nsd = 3.
\end{cases} 
\end{equation}
For $\mathbb{Q}_1$ discretizations, we use $C_{\mathrm{inv}} \approx 1$. To improve the consistency of our formulation in combination with linear finite elements, the second-order derivatives in the residual $\text{res}(u^h)$ are obtained with a least-squares recovery technique~\cite{jansen1999better}.

For the parabolic case ($\ba=\zero$) and linear approximation functions, the weak form in Equation~\eqref{eq:weakAD} is very close to the locally stabilized space-time finite element method presented by Langer and Schafelner in~\cite[Section 3]{langer2019spacetime}. Only the definition of the stabilization parameter, $\tau_\mathrm{SUPG}$ or $\Theta_K$,~\cite[Remark 13.4]{langer2019coeff}, and the enforcement of the initial condition differ.

To provide a first test case for the four space-time discretization methods, we analyze the transient one-dimensional model problem
\begin{align}
\label{eq:IBVP1}
\text{IBVP 1} \quad
\begin{cases}
	\text{res}(u(x,t)) = 0, &\;   x \in ]-1,1[, \quad t \in ]0,2], \\
	u(x,t) = -\sin(\pi x),& \;  \text{at} \; t=0,   \\
	u(-1,t) = u (1,t) , &\;  \text{on} \; P^D. \\
\end{cases}
\end{align}

We consider a time interval $I = \,]0,2]$ and the spatial computational domain $\Omega$ spans from $-1$ to $1$. The model problem is characterized by the periodic boundary conditions and has the analytical solution
\begin{equation}
u(x,t) = -\sin (\pi(x-at)) e^{-k \pi^2 t}.
\end{equation}
The test case setup of IBVP 1 is also discussed by Mojtabi and Deville~\cite{mojtabi2015one} and on a shifted computational domain by Shakib and Hughes~\cite{shakib1991new}. 

\newcommand{\swW}{0.43\textwidth}

\newcommand{\addSWAxis}[2]{
%\tikzsetnextfilename{spacetime-sw#1}
\begin{tikzpicture}[
    axis/.style={very thick, ->},
    every node/.style={color=black}
    ]
       \node at (0,0) [inner sep=0pt, anchor=south west] {\includegraphics[width=\swW,trim={#2cm #2cm 0cm 0cm},clip]{spacetime/SineWave/#1}};
       \draw[axis] (0.0,0.04*\swW)  -- (1.03*\swW,0.04*\swW) node(xline)[right]{$x$};
       \draw[axis] (0.03*\swW,0.0)  -- (0.03*\swW,1.03*\swW) node(xline)[left]{$t$};
\end{tikzpicture}
}

\begin{figure}
\captionsetup[sub]{position=bottom}
\centering 
\includegraphics[width=0.75\textwidth,trim={0cm 0cm 0cm 0cm},clip]{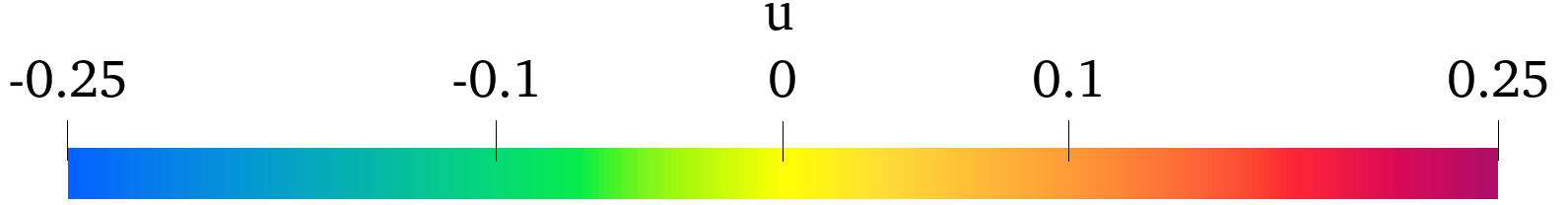}\\ \vspace{1ex}
\subcaptionbox{C-PST.\label{fig:swCPST}}{\addSWAxis{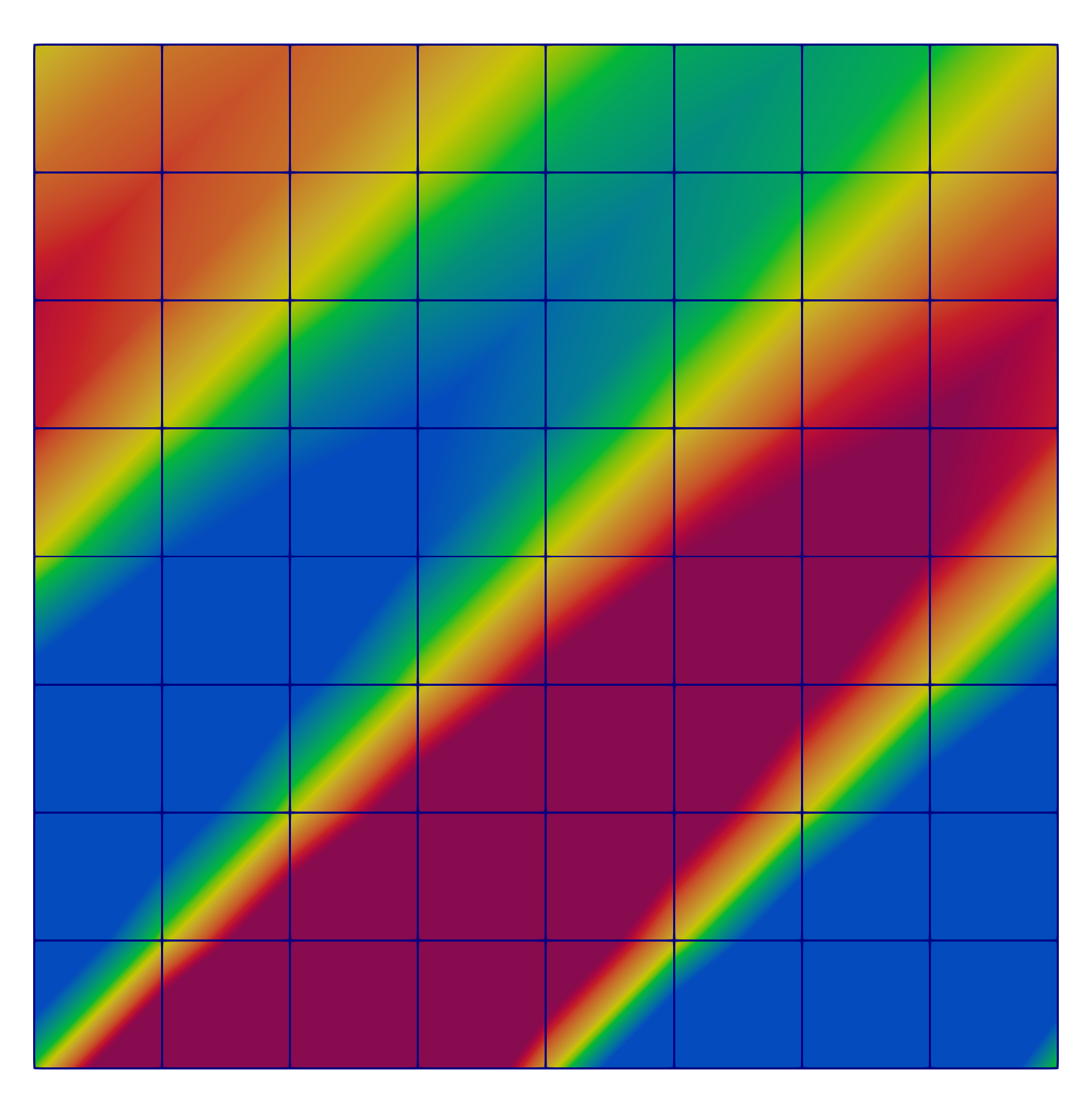}{0}}%
\subcaptionbox{C-SST.\label{fig:swCSST}}{\addSWAxis{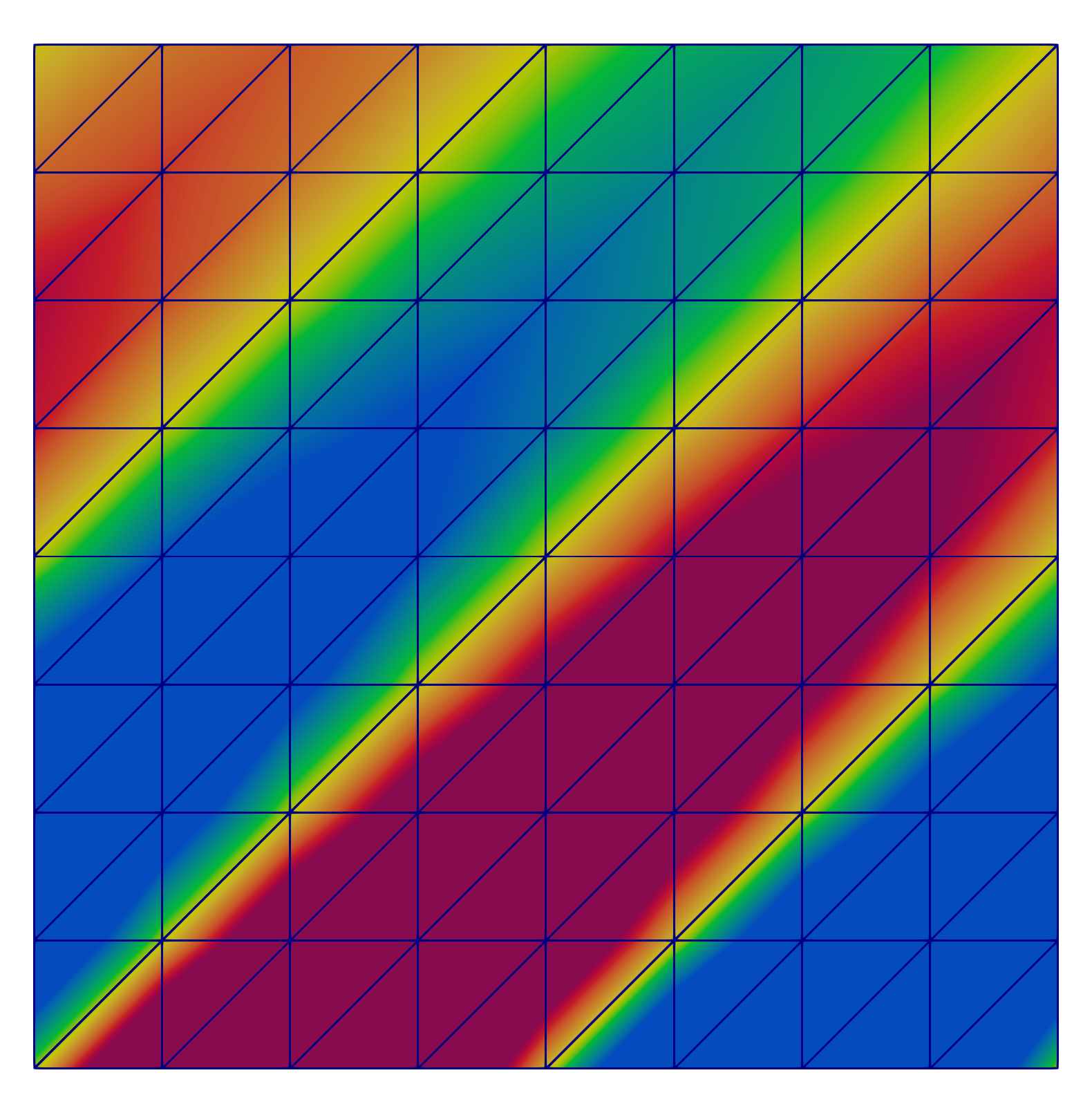}{0}}\\\vspace{1ex}
\subcaptionbox{D-PST.\label{fig:swDPST}}{\addSWAxis{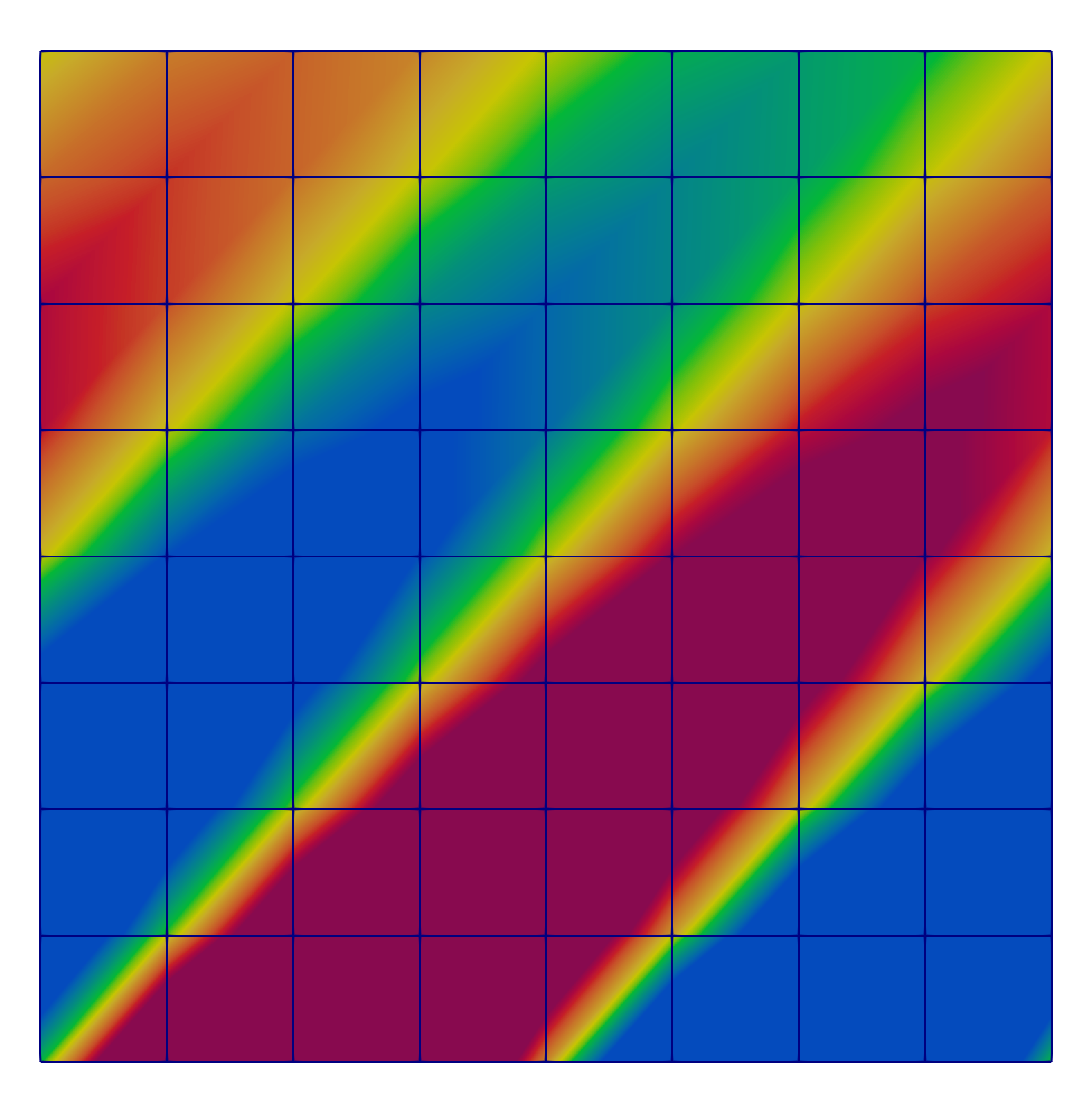}{0.3}}%
\subcaptionbox{D-SST.\label{fig:swDSST}}{\addSWAxis{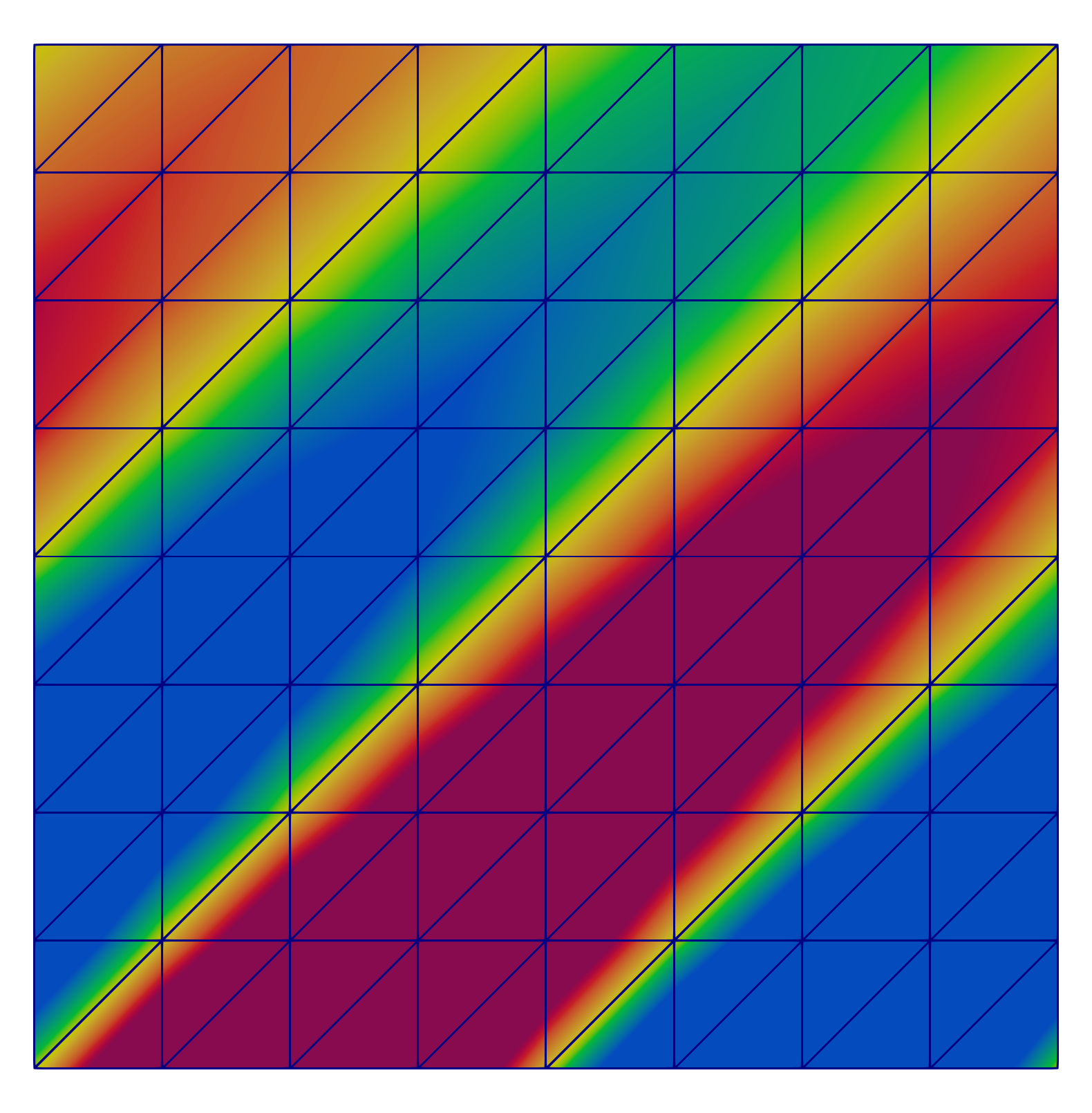}{0}}\\
\caption[Solution of model problem IBVP 1.]{Solution of IBVP 1 computed with four space-time discretization methods.}
\label{fig:sineWave}
\end{figure}

In the numerical solution procedure, we discretize the computational domain $Q$ with eight elements in spatial and temporal direction as shown in Figure~\ref{fig:sineWave}. Due to the periodic boundary conditions  $u(-1,t) = u (1,t)$, this leads to eight independent degrees of freedom in spatial direction. As the initial condition is enforced weakly, the time-continuous discretizations have nine nodes in time direction with one degree of freedom each. Therefore, C-PST and C-SST simulations use $9 \times 8 = 72$ degrees of freedom in total in this specific test case. The time-discontinuous methods have two degrees of freedom per time step to approximate the solution in temporal direction, so $16 \times 8 = 128$ degrees of freedom in total for this specific computation.

An advection speed of $a=1$ and a diffusion coefficient $k=0.1$ lead to the damped traveling sine wave shown in Figure~\ref{fig:sineWave}. For this parameter set, advective and diffusive effects are of similar importance as both are visible to the naked eye. We therefore calibrate the P\'{e}clet number (Equation~\eqref{eq:Peclet}) for this model problem with a characteristic length $L_c = 1/10$ to obtain $\Peclet = 1$ for this parameter set.

Comparing the solution of D-PST in Figure~\ref{fig:swDPST} with the C-PST solution in Figure~\ref{fig:swCPST}, one can note jumps in the solution at the interfaces between the space-time slabs. These small discontinuities in the solution are in line with the weak enforcement of the continuity requirement in the weak form (Equation~\eqref{eq:weakAD}). Also the D-SST solution is discontinuous at the interfaces between space-time slabs. However, these jumps are less pronounced and not visible in the rendering of Figure~\ref{fig:swDSST}. Regarding the SST discretizations (in Figure~\ref{fig:swCSST} and~\ref{fig:swDSST}), we can note that the solution $u$ is advected along the diagonal edges of the SST discretizations. In this particular case with $\Delta x= \Delta t$ and $a=1$, the characteristics perfectly align with the finite element edges.

\newcommand{\plthgtC}{0.2*\textheight}
\newcommand{\pltwdtC}{0.37*\textwidth}
\renewcommand{\pltwdtC}{7cm}

\tikzset{dataline/.style={no markers, very thick}}

\newcommand{\PlotDifference}{
%\tikzsetnextfilename{spacetime-difference}
\begin{tikzpicture}
\begin{axis}[
     	height=\plthgtC, width=\pltwdtC,
			scale only axis,
	xmin = -1,
	xmax = 1,
	ymin = -0.027,
	ymax = 0.027,
	xlabel={$x$},
   	ylabel={$(u-u^h)(\cdot,t_f)$},
	y label style={at={(axis description cs:-0.075,.5)}},
	y tick scale label style={xshift=-0.5cm},
	grid=major,
       	cycle list name=myQPlotCycleList,
]
		\addplot+[dataline] table [x index=4, y expr = -sin(180*\thisrowno{4})*0.1389111331-\thisrowno{1}, each nth point=1, col sep=comma] {spacetime/cpst.csv};
		\addplot+[dataline] table [x index=5, y expr = -sin(180*\thisrowno{5})*0.1389111331-\thisrowno{2}, each nth point=1, col sep=comma] {spacetime/csst.csv};
	        \addplot+[dataline] table [x index=4, y expr = -sin(180*\thisrowno{4})*0.1389111331-\thisrowno{1}, each nth point=1, col sep=comma] {spacetime/dpst.csv};
	        \addplot+[dataline] table [x index=4, y expr = -sin(180*\thisrowno{4})*0.1389111331-\thisrowno{1}, each nth point=1, col sep=comma] {spacetime/dsst.csv};
\end{axis}
\end{tikzpicture}
}

\newcommand{\PlotDifferenceAtNodes}{
%\tikzsetnextfilename{spacetime-differenceAtNodes}
\begin{tikzpicture}
\begin{axis}[
       height=\plthgtC, width=\pltwdtC,
       		scale only axis,
		xmin = -1,
		xmax = 1,
		ymin = -0.027,
		ymax = 0.027,
		xlabel={$x$},
   		ylabel={$(u-u^h)(x^h,t_f)$},
		y label style={at={(axis description cs:-0.075,.5)}},
		y tick scale label style={xshift=-0.5cm},
		grid=major,
		cycle list name=myQPlotCycleList,
               legend style={
               at={(0.5,1.3)},
               anchor=north,
               cells={anchor=west},
               legend columns=4,
       },
       legend entries={C-PST, C-SST, D-PST, D-SST},
       legend to name=differenceLegend,
]
	\addplot+[dataline] table [x index=4, y expr = \thisrowno{0}, each nth point=1, col sep=comma] {spacetime/cpst.csv};
	\addplot+[dataline] table [x index=5, y expr = \thisrowno{0}, each nth point=1, col sep=comma] {spacetime/csst.csv};
	\addplot+[dataline] table [x index=4, y expr = \thisrowno{0}, each nth point=1, col sep=comma] {spacetime/dpst.csv};
	\addplot+[dataline] table [x index=4, y expr = \thisrowno{0}, each nth point=1, col sep=comma] {spacetime/dsst.csv};
\end{axis}
\end{tikzpicture}
}

\begin{figure}
\centering
\parbox{.49\linewidth}{%
\flushright \parbox{\pltwdtC}{\centering \subcaption{Differences between $u$ and $u^h$.\label{fig:diffFunction}}}%
\hfill \PlotDifference}\,\parbox{.49\linewidth}{\flushright \parbox{\pltwdtC}{\centering \subcaption{Differences at nodes $x^h$.\label{fig:diffNodes}}}%
\hfill \PlotDifferenceAtNodes}\\
%\ref{differenceLegend}\\
\includegraphics{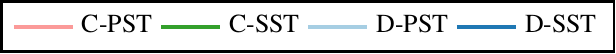}
\caption{Comparison of four space-time discretization methods for IBVP 1.}
\label{fig:difference}
\end{figure}

Figure~\ref{fig:difference} compares the numerical solutions $u^h$ of the four space-time discretization methods to the analytical solution $u$ at the final time $t_f$. In the plot of the differences $u-u^h$ (Figure~\ref{fig:diffFunction}), the interpolation error between the nodal values is very prominent. Please, note that this error is inherent to the linear interpolation of a trigonometric function. Removing this unavoidable error (for linear approximation functions), Figure~\ref{fig:diffNodes} connects the values at the finite element nodes with straight line segments. For the employed, very coarse discretizations, the nodal differences of the SST solutions to the analytical solution are smaller, despite the smaller number of degrees of freedom in comparison to the PST methods. Additionally, the D-SST method shows hardly any phase error. Comparing Figure~\ref{fig:diffFunction} and Figure~\ref{fig:diffNodes}, one can observe that the error of the finite element solution at the nodes is of the same order as the interpolation error. 

Returning to the complete space-time solution (Figure~\ref{fig:sineWave}), all four space-time discretization methods arrive at similar results. Given the extremely coarse discretization, one can consider all numerical solutions to be in accordance with the analytical solution. We therefore conclude that all four space-time discretization schemes are suitable for advection-diffusion problems.

The C-PST method has been analyzed for the heat equation theoretically and with numerical experiments by Aziz and Monk~\cite{aziz1989continuous}. It is found that the use of linear finite element approximation functions in C-PST leads to a version of the Crank-Nicolson method. Moreover, the tensor-product approach of C-PST leads to a global linear equation system with specific structure for parabolic initial boundary value problems. This structure can be exploited in the construction of an efficient parallel solver as shown by Langer and Zank~\cite{langer2020efficient}. Still, we will not further consider the scheme in this paper.

The C-SST method allows for space-time adaptivity on unstructured meshes~\cite{langer2019spacetime} and in Section~\ref{sec:prp} the C-SST method is used to include topology changes of the spatial computational domain $\Omega$ in a boundary-conforming space-time mesh. However, Section~\ref{sec:convergence} focuses on the time-discontinuous methods, D-PST and D-SST.

\section{Computational Error Analysis of Time-Discontinuous Discretizations}
\label{sec:convergence}

To investigate the convergence of the time-discontinuous space-time discretizations D-PST and D-SST, a computational error analysis is performed. The following space-time convergence studies consider two test cases. Before investigating IBVP 1 (Equation~\eqref{eq:IBVP1}) for six model parameter sets in Section~\ref{ssec:COne}, we first consider the parabolic case ($a = 0$) of a second initial boundary value problem IBVP 2 with time-dependent Dirichlet boundary conditions in Section~\ref{ssec:CTwo}. Both initial boundary value problems have analytical solutions, which serve particularly well as reference solutions in the convergence studies, since they are independent of implementation issues or round-off errors introduced in computer arithmetic.

For each model problem, parameter set, and discretization method (D-PST, D-SST) a space-time convergence study with 198 simulations is performed. The numerical simulation settings are obtained as follows. We divide the computational domain in time direction (up to the final time $t_f = 2$) into $\nts$ space-time slabs of constant size $\Delta t = t_f / \nts$. We consider 15 levels of recursive temporal refinement such that $\nts$ is doubled from the coarser to the finer level 
\begin{equation} 
\nts = 2^{m-1},\, m = 4, \dots, 18.
\end{equation}
In the same manner, the spatial domain is divided into $\nex$ elements of constant size $\Delta x = 2 / \nex$. The number of elements in spatial direction is given by
\begin{equation}
\nex = 2^{l-1}, \,l = 4, \dots, 18.
\end{equation}

For each simulation, the relative $L_2$ error $e_{lm}$ at the final time $t_f=2$ for the spatial refinement level $l$ and temporal refinement level $m$ is evaluated. In practice, we use an element-wise two-point Gaussian quadrature for the spatial integration
\begin{align}
\label{eq:L2ad}
& e_{lm} = \frac{\| (u - u^h)(\cdot,t_f) \|}{\| u(\cdot,t_f) \|} \approx \frac{1}{\| u(\cdot,t_f) \|}  \sqrt{ \frac{\Delta x}{2} \sum_{e=1}^{\nex}  \sum_{iq=1}^{2} \left( u(x^e_{iq},t_f) - u^h(x^e_{iq},t_f) \right)^2 }.
\end{align}
 Additionally, we measure the nodal errors as 
 \begin{align}
 \label{eq:adNe}
& E_{lm} \coloneqq \frac{1}{\| u(\cdot,t_f) \|} \sqrt{ \Delta x \sum_{i=1}^{\nex} \left( u(x_i,t_f) - u^h(x_i,t_f) \right)^2 },
 \end{align}
with the index $i$ running over all nodes, except for the last one. In model problem IBVP 1, node 1 and node $\nex +1$ have identical solution values enforced by the periodic boundary conditions. For both model problems, the prefactor with the $L_2$ norm of the solution at the final time $t_f = 2$ reads \begin{equation} \frac{1}{\| u(\cdot,t_f) \|} = e^{2 k \pi^2}. \end{equation}

\newcommand{\plthgtTDT}{0.33\textheight}
\newcommand{\pltwdtTDT}{0.52\textwidth}
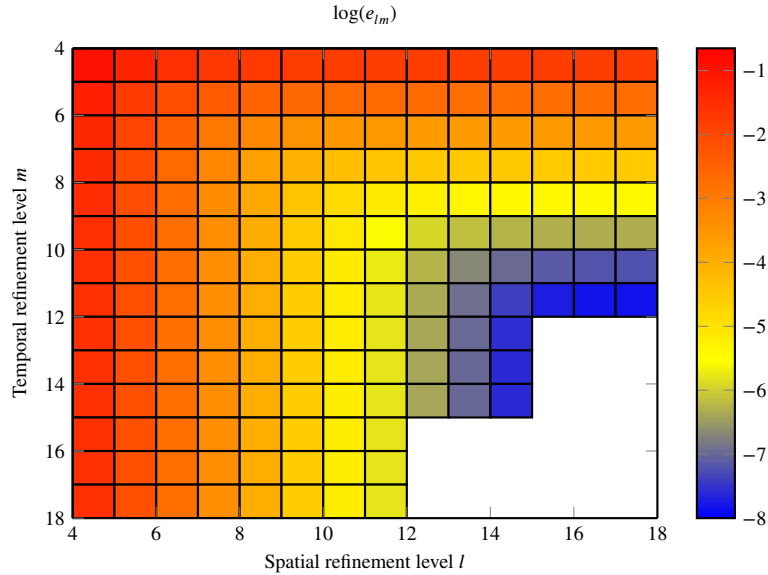
\begin{figure}
\centering
{\footnotesize
%\tikzsetnextfilename{spacetime-orga}
\begin{tikzpicture}
\begin{axis}
[
height=\plthgtTDT, width=\pltwdtTDT,
point meta min = -8,
zmax=0,
zmin=-12,
unbounded coords=jump,
view={0}{90},
y dir=reverse,
colorbar right,
draw=black,
%colorbar sampled,
%colorbar style={samples=25},
xlabel=Spatial refinement level $l$, ylabel=Temporal refinement level $m$, 
title=$\log (e_{lm})$,
]
%\addplot3 [surf, shader=flat, draw=black, mesh/rows=15, mesh/ordering=y varies, z filter/.expression= {z==-15 ? nan : z}] table [col sep=comma]  {Data/FSTQuadNu0.0L2ErrorCoords.csv};
\addplot3 [surf, shader=flat, draw=black, mesh/rows=15, mesh/ordering=y varies, z filter/.expression= {z==-15 ? nan : z}] table [col sep=comma]  {Data/FSTQuadNu0.0RealL2ErrorCoords.csv};
%\addplot3 [surf, mesh/rows=15, mesh/ordering=y varies, z filter/.expression= {z==-15 ? nan : z}] table [col sep=comma]  {Data/FSTQuadNu0.0L2ErrorCoords.csv};
\end{axis}
\end{tikzpicture}
}
\caption{Organization of space-time convergence study based on $L_2$ error $e_{lm}$.}
\label{fig:orga}
\end{figure}

Simulations are performed for whole-numbered parameter pairs $(l,m)$ corresponding to grid line intersections in Figure~\ref{fig:orga}. To avoid unnecessary computational cost, we omit combinations of the finest refinement levels as shown in Figure~\ref{fig:orga}. Note that the patch color is based on the mean value of the $L_2$ error $e_{lm}$ of the four simulations connected by a patch. 

To check for spatial convergence, we consider the finest temporal refinement level $m= 18$ and vary $l = 4, \dots, 12$, which corresponds to the bottom line of the plot in Figure~\ref{fig:orga}. Analogously to investigate temporal convergence, we consider the finest spatial refinement level $l= 18$ and vary $m = 4, \dots, 12$. This corresponds to the rightmost line of the plot in Figure~\ref{fig:orga}. On the space-time diagonal $l=m$, the numerical values of $\Delta t$ and $\Delta x$ coincide. Despite the different units that one would assign to the physical quantities, we use $\Delta t = \Delta x$ to express that the numerical values are equal. Along the curve $\Delta t= \Delta x$, twelve data points $m =l = 4, \dots, 15$ are generated.

%%% IBVP 2
\subsection{Parabolic model problem IBVP 2}\label{ssec:CTwo}

In this section, we study the pure diffusion case of the model problem
\begin{align}
\label{eq:IBVP2}
\text{IBVP 2} \quad
\begin{cases}
	\text{res}(u(x,t)) = 0, &\;   x \in ]-1,1[, \quad t \in ]0,2], \\
	u(x,t) = \cos(\pi x),& \;  \text{at} \; t=0,   \\
	u(-1,t) = u (1,t) = b(t) = - e^{-k \pi^2 t}, &\;  \text{on} \; P^D. \\
\end{cases}
\end{align}
With the time-dependent Dirichlet boundary conditions $b(t)$, IBVP 2 has the analytical solution
\begin{equation}
u(x,t) = \cos (\pi x) e^{-k \pi^2 t}.  \label{eq:P2sol}
\end{equation}

The considered advection-diffusion equation~\eqref{eq:ADstrong} reduces for $a=0$ to the heat equation. For a corresponding initial boundary value problem with homogeneous boundary conditions, convergence estimates for D-PST discretizations are known from literature. Thom\'{e}e presents in~\cite[Theorem 12.7 ]{thomee2006Galerkin} a superconvergence result for the temporal discretization error at the final time $t_f$. Considering linear basis functions, the error bound for the parabolic problem can be summarized as
\begin{equation}
\label{eq:thomee}
\| (u - u^h)(\cdot,t_f) \| \leq C (\Delta t^3 + \Delta x^2),
\end{equation}
where $C$ is a positive constant independent of $\Delta t$ and $\Delta x$.

\newcommand{\plthgtTD}{0.33\textheight}
\newcommand{\pltwdtTD}{0.5\textwidth}

\pgfplotsset{
        small e axis style/.style={
       		ymin = 2e-10, ymax = 0.5,
       		ylabel= {$e_{lm}$},
		y label style={at={(axis description cs:0.0,1.0)},,rotate=-90,anchor=south east},
		grid=major,
		scale only axis,
        },
        large e axis style/.style={
       		ymin = 2e-10, ymax = 0.5,
       		ylabel= {$E_{lm}$},
		grid=major,
        }
}

\newcommand{\logPlotThreeD}[1]{
{\footnotesize
%\tikzsetnextfilename{spacetime-#1}
\begin{tikzpicture}
\begin{axis}
[
height=\plthgtTD, width=\pltwdtTD,
point meta min = -8,
zmax=0,
zmin=-12,
unbounded coords=jump,
y dir=reverse,
xlabel=$l$, ylabel=$m$, 
zlabel=$\log (e_{lm})$,
z label style={at={(axis description cs:0.0,0.95)},,rotate=-90,anchor=north},
]
\addplot3 [surf, mesh/rows=15, mesh/ordering=y varies, z filter/.expression= {z==-15 ? nan : z}] table [col sep=comma]  {Data/#1RealL2ErrorCoords.csv};
\end{axis}
\end{tikzpicture}
}
}

\begin{figure}
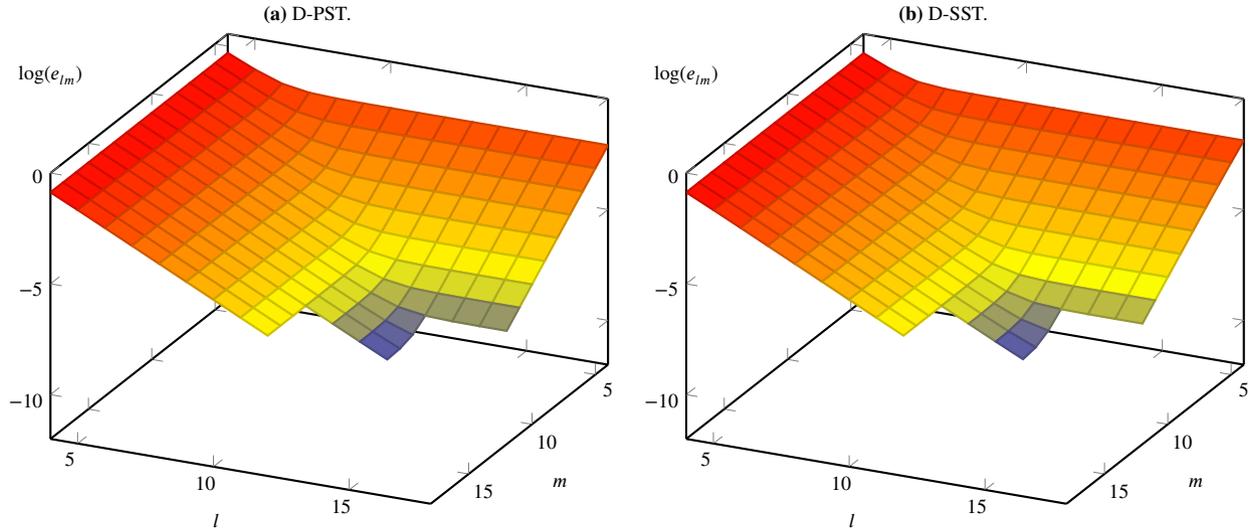

\centering
\subcaptionbox{D-PST. \label{fig:convFSTP2}}{\logPlotThreeD{FSTQuadNu0.1P2}}%
\subcaptionbox{D-SST. \label{fig:convSSTP2}}{\logPlotThreeD{SST1DNu0.1P2}} \\
\caption[Convergence of $L_2$ error for IBVP 2.]{Convergence visualization of $L_2$ error for parabolic problem configuration IBVP 2.}
\label{fig:convCTwo}
\end{figure} 

In the following, we compare our computational findings to the theoretical result above. The results of the space-time convergence studies are visualized in convergence surfaces (see Figure~\ref{fig:convCTwo}). The surfaces are obtained by plotting the $L_2$ error $e_{lm}$ in logarithmic scale over the spatial and temporal refinement level indices $l$ and $m$. Corresponding convergence surfaces based on the nodal error measure can be found in the Appendix~\ref{sec:appendix} in Figure~\ref{fig:appC2Ne}. For both discretization methods, the error plots result in a continuous surface (Figure~\ref{fig:convFSTP2} and~\ref{fig:convSSTP2}). Moreover, the surfaces show, that the error values in the area of the diagonal ($l=m$, $\Delta t = \Delta x$) are influenced by the spatial and temporal mesh size. However, on the finest spatial discretization level ($l=18$), the error varies only with $\Delta t$. The same holds for the finest temporal refinement level ($m=18$) and $\Delta x$. Therefore, extracting the curves $l=18$ or $m=18$ from the convergence surfaces gives us the isolated spatial or temporal convergence behavior of the methods.

\newcommand{\logPlotTwoDPTwoSpace}{
\tikzsetnextfilename{spacetime-TwoDPTwoSpace}
\begin{tikzpicture}
\begin{loglogaxis}[
      height=\plthgtC, width=\pltwdtC,
       xlabel={$\Delta x$},
       small e axis style,
       		cycle list name=myQPlotCycleList,
]
	\addplot+[dataline] table [x expr=2/(2^(\coordindex+3)), y index=0] {Data/FSTQuadNu0.1P2RealL2ErrorSpace.txt};
	\addplot+[dataline] table [x expr=2/(2^(\coordindex+3)), y index=0] {Data/SST1DNu0.1P2RealL2ErrorSpace.txt};
	\addplot[very thick, dashed, domain=0.001:0.25, samples=100]{(x)^(2)};
	\addplot+[dataline] table [x expr=2/(2^(\coordindex+3)), y index=0] {Data/FSTQuadNu0.1P2ModRealL2ErrorSpace.txt};
	\addplot+[dataline] table [x expr=2/(2^(\coordindex+3)), y index=0] {Data/SST1DNu0.1P2ModRealL2ErrorSpace.txt};
\end{loglogaxis}
\end{tikzpicture}
}

\newcommand{\logPlotTwoDPTwoTime}{
%\tikzsetnextfilename{spacetime-TwoDPTwoTime}
\begin{tikzpicture}
\begin{loglogaxis}[
       	height=\plthgtC, width=\pltwdtC,
       	xlabel={$ \Delta t$},
	small e axis style,
	cycle list name=myQPlotCycleList,
       	legend style={
               at={(0.5,1.3)},
               anchor=north,
               cells={anchor=west},
               legend columns=3,
       },
       legend entries={D-PST,  D-SST, $\mathcal{O}(\Delta x^2 )$ or $\mathcal{O}(\Delta t^2 )$,  D-PST with $\tilde{b}(t_l)$,  D-SST with $\tilde{b}(t_l)$ ,$\mathcal{O}(\Delta t^3)$},
       legend to name=twoDPTwoLegend,
]
	\addplot+[dataline] table [x expr=2/(2^(\coordindex+3)), y index=0] {Data/FSTQuadNu0.1P2RealL2ErrorTime.txt};
	\addplot+[dataline] table [x expr=2/(2^(\coordindex+3)), y index=0] {Data/SST1DNu0.1P2RealL2ErrorTime.txt};
	\addplot[very thick, dashed, domain=0.001:0.25, samples=100]{(x)^(2)/15};
	\addplot+[dataline] table [x expr=2/(2^(\coordindex+3)), y index=0] {Data/FSTQuadNu0.1P2ModRealL2ErrorTime.txt};
	\addplot+[dataline] table [x expr=2/(2^(\coordindex+3)), y index=0] {Data/SST1DNu0.1P2ModRealL2ErrorTime.txt};
         \addplot[very thick, dotted, domain=0.01:0.25, samples=100]{(x)^(3)/100};
\end{loglogaxis}
\end{tikzpicture}
}

%\begin{figure}
%\centering
%\ref{twoDPTwoLegend}\\
%\subcaptionbox{Spatial convergence $\Delta t = \frac{1}{65536}$ \label{fig:convTwoDPTwoSpace}}{\logPlotTwoDPTwoSpace}\,
%\subcaptionbox{Temporal convergence $\Delta x = \frac{1}{65536}$ \label{fig:convTwoDPTwoTime}}{\logPlotTwoDPTwoTime}
%\caption{Spatial and temporal convergence for IBVP 2.}
%\label{fig:twoDPTwo}
%\end{figure} 

\begin{figure}
\centering
\parbox{.49\linewidth}{
\flushright \parbox{\pltwdtC}{\centering \subcaption{Spatial convergence, $\Delta t = \frac{1}{65536}$. \label{fig:convTwoDPTwoSpace} }}
\hfill \logPlotTwoDPTwoSpace
}\,\parbox{.49\linewidth}{
\flushright \parbox{\pltwdtC}{\centering \subcaption{Temporal convergence, $\Delta x = \frac{1}{65536}$. \label{fig:convTwoDPTwoTime}}}
\hfill \logPlotTwoDPTwoTime}
%\ref{twoDPTwoLegend}\\
\includegraphics{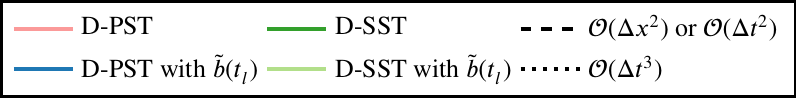}
\caption{Spatial and temporal convergence for IBVP 2.}
\label{fig:twoDPTwo}
\end{figure} 

At first, focusing on the spatial convergence rates shown in Figure~\ref{fig:convTwoDPTwoSpace}, we observe a second-order spatial accuracy for both methods as the curves of D-PST and D-SST coincide. This is to be expected, as the same spatial mesh is used. Moreover, this observation is also in line with the theoretical result given in Equation~\eqref{eq:thomee}. Next, looking at the temporal convergence rates in Figure~\ref{fig:convTwoDPTwoTime}, a second order temporal convergence is observed for both methods. This is in strict contrast to the third-order time accuracy expected from Equation~\eqref{eq:thomee}.

As pointed out by Voulis and Reusken in~\cite{voulis2019discontinuous}, the reduced convergence order is due to the time-dependent boundary conditions. Moreover, it is shown in their work that superconvergence can be recovered by applying a temporal interpolation operator to the boundary condition $b(t)$. The use of this interpolated boundary condition is equivalent to the time-discontinuous discretization of the boundary condition $\frac{\partial u}{\partial t}(\bx,t)=\frac{\partial b}{\partial t}(t)$.  In our considered test case, the temporal convergence can be improved with the following treatment. On the upper time level $\Omega_u=\Omega(t_u)$ of each space-time slab, the boundary condition is precisely evaluated as
\begin{equation}
u(-1,t_u) = u (1,t_u) = b(t_u) = - e^{-k \pi^2 t_u}.
\end{equation}
On the lower time level $\Omega_l$, a modified boundary condition $\tilde{b}(t_l)$ is applied. The modified boundary condition is constructed such that the linear interpolation of the finite element shape functions leads to the correct analytical mean of the boundary condition
\begin{equation}
\frac{1}{2} \left[ b(t_u) + \tilde{b}(t_l) \right] \stackrel{!}{=} \frac{1}{\Delta t} \int_{t_l}^{t_u} b(t) dt.
\end{equation}
For the considered example, this yields
\begin{equation}
\tilde{b}(t_l) = \left[1 + \frac{2}{k \pi^{2} \Delta t } \left( 1 -  e^{k \pi^{2} \Delta t} \right)\right] e^{-k \pi^{2} t_u}.
\end{equation}

\begin{figure}
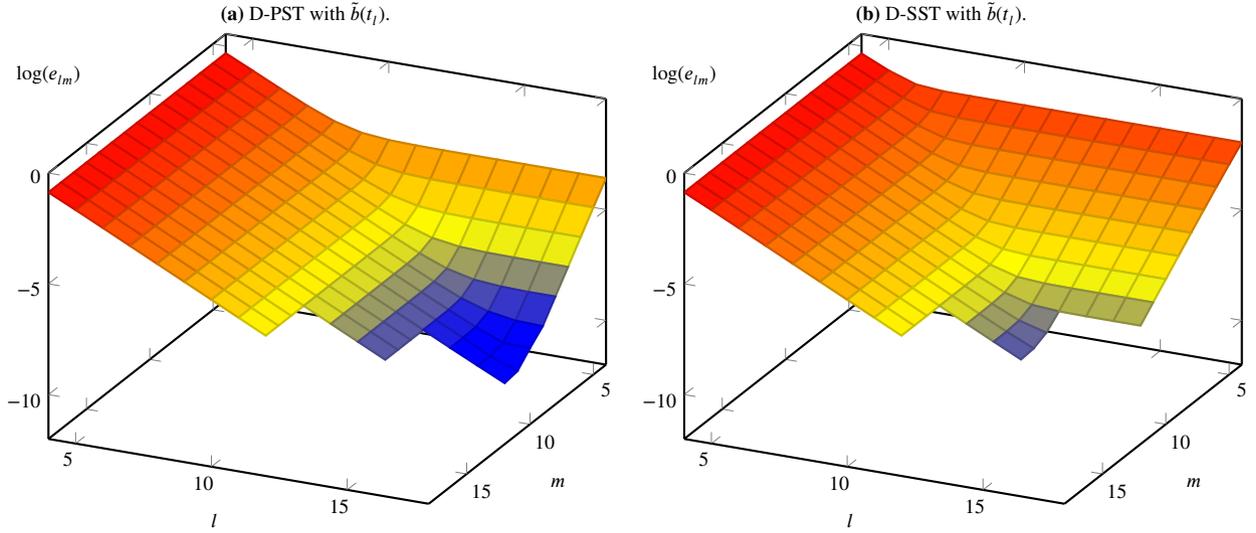

\centering
\subcaptionbox{D-PST with $\tilde{b}(t_l)$. \label{fig:convFSTP2Mod}}{\logPlotThreeD{FSTQuadNu0.1P2Mod}}%
\subcaptionbox{D-SST with $\tilde{b}(t_l)$. \label{fig:convSSTP2Mod}}{\logPlotThreeD{SST1DNu0.1P2Mod}} \\
\caption[Convergence of $L_2$ error for IBVP 2 with $\tilde{b}(t_l)$.]{Convergence visualization of $L_2$ error for IBVP 2 with modified boundary condition.}
\label{fig:convCTwoTilde}
\end{figure} 

Repeating the space-time convergence study with modified boundary conditions, we obtain the results shown in Figure~\ref{fig:convCTwoTilde}. Here, D-PST reaches significantly smaller error values in comparison to the case shown in Figure~\ref{fig:convFSTP2}. Returning in the line plots of Figure~\ref{fig:convTwoDPTwoSpace}, it can be seen that the second order spatial convergence of both methods is not affected by the boundary condition treatment as all four curves coincide. But, for D-PST with $\tilde{b}(t_l)$, third-order temporal convergence is indeed obtained (Figure~\ref{fig:convTwoDPTwoTime}). The result numerically confirms that temporal superconvergence (as stated in Equation~\eqref{eq:thomee}) can also be obtained for time-dependent boundary conditions with a proper treatment~\cite{voulis2019discontinuous}.
 
For D-SST with $\tilde{b}(t_l)$, only quadratic temporal convergence is observed. The lower convergence order of the D-SST method with treatment of the time-dependent boundary conditions hints at the fact that superconvergence of the D-PST method is linked to the tensor-product structure of the discretization. However, also in case of the D-SST discretization, the proposed treatment of time-dependent boundary conditions is helpful\textemdash the error values decrease by approximately 25\%. 

The purely spatial or temporal refinements are interesting as they show an isolated spatial or temporal convergence behavior, but they are certainly not efficient in terms of computational cost that is required to obtain a certain level of accuracy. Elaborating on this, we assume that the computational cost of a simulation is related to the number of degrees of freedom $\ndf$. For the considered discretizations, $\ndf$ can be expressed by the number of time steps $\nts$ and the number or elements in $x$-direction $\nex$ as
\begin{equation}
\ndf = 2 \cdot \nts \cdot \left(\nex -1 \right).
\end{equation}

Since we estimate the total computational cost by the number of degrees of freedom, it is of the order $\mathcal{O}(\frac{1}{\Delta x}\frac{1}{\Delta t})$. This can be used to balance the spatial and temporal discretization to minimize the computational cost for a desired error. The optimal relation between the spatial mesh size and the temporal therefore depends on the relation between the spatial and the temporal convergence order. If the spatial and temporal convergence order match, then the choice $\Delta x = \Delta t$ is optimal.  However, if we consider the setting in Equation \eqref{eq:thomee}, then the optimal choice is $\Delta x^2 = \Delta t^3$.

In the visualizations of the convergence studies, e.g., in Figure~\ref{fig:convCTwoTilde}, we can identify the best space-time refinement strategy as the steepest decent in the convergence surfaces. In Figure~\ref{fig:convSSTP2Mod}, an advantageous space-time refinement strategy for D-SST essentially follows $\Delta t = 4 \cdot \Delta x$. In contrast, for D-PST with $\tilde{b}(t_l)$ (Figure~\ref{fig:convFSTP2Mod}), the second-order spatial accuracy and third-order temporal accuracy lead to an advantageous space-time refinement strategy along the curve $\Delta t^3 = \Delta t^2$. 

\newcommand{\logPlotTwoDPTwoSpaceTime}{
%\tikzsetnextfilename{spacetime-TwoDPTwoSpaceTime}
\begin{tikzpicture}
\begin{loglogaxis}[
       	height=\plthgtC, width=\pltwdtC,
       	xlabel={$ \Delta t$},
	small e axis style,
		cycle list name=myQPlotCycleList,
]
	\addplot+[dataline] table [x expr=2/(2^(\coordindex+5)), y index=0] {Data/FSTQuadNu0.1P2RealL2Error.csvOffset2};
	\addplot+[dataline] table [x expr=2/(2^(\coordindex+5)), y index=0] {Data/SST1DNu0.1P2RealL2Error.csvOffset2};	
	\addplot[very thick, dashed, domain=0.000125:0.065, samples=100]{(x)^(2)/5};
	\addplot+[dataline] table [x expr=2/(2^(\coordindex+5)), y index=0] {Data/FSTQuadNu0.1P2ModRealL2Error.csvOffset2};
	\addplot+[dataline] table [x expr=2/(2^(\coordindex+5)), y index=0] {Data/SST1DNu0.1P2ModRealL2Error.csvOffset2};

\end{loglogaxis}
\end{tikzpicture}
}

\newcommand{\logPlotTwoDPTwoSpaceTimeSpecial}{
%\tikzsetnextfilename{spacetime-TwoDPTwoSpaceTimeSpecial}
\begin{tikzpicture}
\begin{loglogaxis}[
       	height=\plthgtC, width=\pltwdtC,
       	xlabel={$ \Delta t$},
	small e axis style,
       	grid=major,
	cycle list name=myQPlotCycleList,
	legend style={
               at={(0.5,1.3)},
               anchor=north,
               cells={anchor=west},
               legend columns=3,
       },
       legend entries={D-PST, D-SST, $\mathcal{O}(\Delta t^2 )$, D-PST with $\tilde{b}(t_l)$, D-SST with $\tilde{b}(t_l)$, $\mathcal{O}(\Delta t^3)$},
       legend to name=twoDPTwoLegendST,
]
	\addplot+[dataline] table [x expr=2^(-2*\coordindex -2), y index=0] {Data/FSTQuadNu0.1P2RealL2ErrorSpaceTimeSpecial.txt};
	\addplot+[dataline] table [x expr=2^(-2*\coordindex -2), y index=0] {Data/SST1DNu0.1P2RealL2ErrorSpaceTimeSpecial.txt};
	\addplot[very thick, dashed, domain=0.001:0.025, samples=100]{(x)^(2)/10};
	\addplot+[dataline] table [x expr=2^(-2*\coordindex -2), y index=0] {Data/FSTQuadNu0.1P2ModRealL2ErrorSpaceTimeSpecial.txt};
	\addplot+[dataline] table [x expr=2^(-2*\coordindex -2), y index=0] {Data/SST1DNu0.1P2ModRealL2ErrorSpaceTimeSpecial.txt};
	 \addplot[very thick, dotted, domain=0.001:0.25, samples=100]{(x)^(3)};
\end{loglogaxis}
\end{tikzpicture}
}

\begin{figure}
\centering
\parbox{.49\linewidth}{
\flushright \parbox{\pltwdtC}{\centering \subcaption{Refinement along $\Delta t = 4 \cdot \Delta x$. \label{fig:convTwoDPTwoSpaceTime} }}
\hfill \logPlotTwoDPTwoSpaceTime
}\,\parbox{.49\linewidth}{
\flushright \parbox{\pltwdtC}{\centering \subcaption{Refinement along $\Delta t^3 = \Delta x^2$. \label{fig:convTwoDPTwoSpaceTimeSpecial}}}
\hfill\logPlotTwoDPTwoSpaceTimeSpecial}
%\ref{twoDPTwoLegendST}\\
\includegraphics{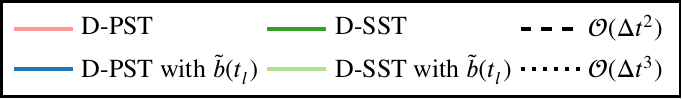}
\caption{Influence of $\tilde{b}(t_l)$ on space-time convergence for IBVP 2.}
\label{fig:twoDPTwoST}
\end{figure} 

The $L_2$ errors along the curves $\Delta t= 4 \cdot \Delta x$ and $\Delta t^3 = \Delta x^2$ are extracted from the convergence surfaces and plotted in Figure~\ref{fig:twoDPTwoST}. On the space-time diagonal with offset ($\Delta t= 4 \cdot \Delta x$, Figure~\ref{fig:convTwoDPTwoSpaceTime}), both methods show a second-order convergence for the computations with and without $\tilde{b}(t_l)$. For D-PST with $\tilde{b}(t_l)$, the curve lies is the zone where the spatial error dominates, hence, we expect second-order convergence also for this method. However, in contrast to Figure~\ref{fig:convTwoDPTwoSpace}, the curves do not coincide and the treatment of the time-dependent boundary conditions proves advantageous in terms of the absolute error values. Note that transitioning from one data point to the next along the space-time diagonal doubles $\nts$ and $\nex$.

Following the advantageous refinement strategy for D-PST with $\tilde{b}(t_l)$, Figure~\ref{fig:convTwoDPTwoSpaceTimeSpecial} shows that the third-order temporal accuracy of D-PST with $\tilde{b}(t_l)$ is retained along the curve $\Delta t^3 = \Delta x^2$. Summarizing the parabolic model problem analysis, both methods, D-PST and D-SST, converge at least quadratically against the analytical solution. With proper treatment, D-PST converges cubically with respect to $\Delta t$ even for time-dependent boundary conditions. 

%%%% Returning to IBVP 1
\subsection{Advective-diffusive model problem IBVP 1}\label{ssec:COne}

While the main challenge in the previous Section~\ref{ssec:CTwo} was the treatment of time-dependent boundary conditions, this section investigates the convergence behavior of the methods as the model parameters transition from the parabolic case to advection-diffusion cases and to the hyperbolic case. The numerical error analysis of D-PST and D-SST is therefore continued with the model problem IBVP 1 (Equation~\ref{eq:IBVP1}). Six parameter sets are considered. They include the parabolic case $a=0$, $k=0.1$, $\Peclet = 0$, next to four advection-diffusion cases with decreasing viscosity $ a=1$, $k = 0.1, 0.01, 0.001, 0.0001$, $\Peclet = 1, 10, 100, 1000$ and the hyperbolic case $a = 1$, $k = 0$, $\Peclet = \infty$. The periodic boundary conditions do not require the treatment of time-dependent boundary conditions. %Note the value of the solution $u(-1,t)$ changes over time, but the condition $u(-1,t)=u(1,t)$ remains constant. ???Do we want to include an explanation? => ASK IGOR??? 

\begin{figure}
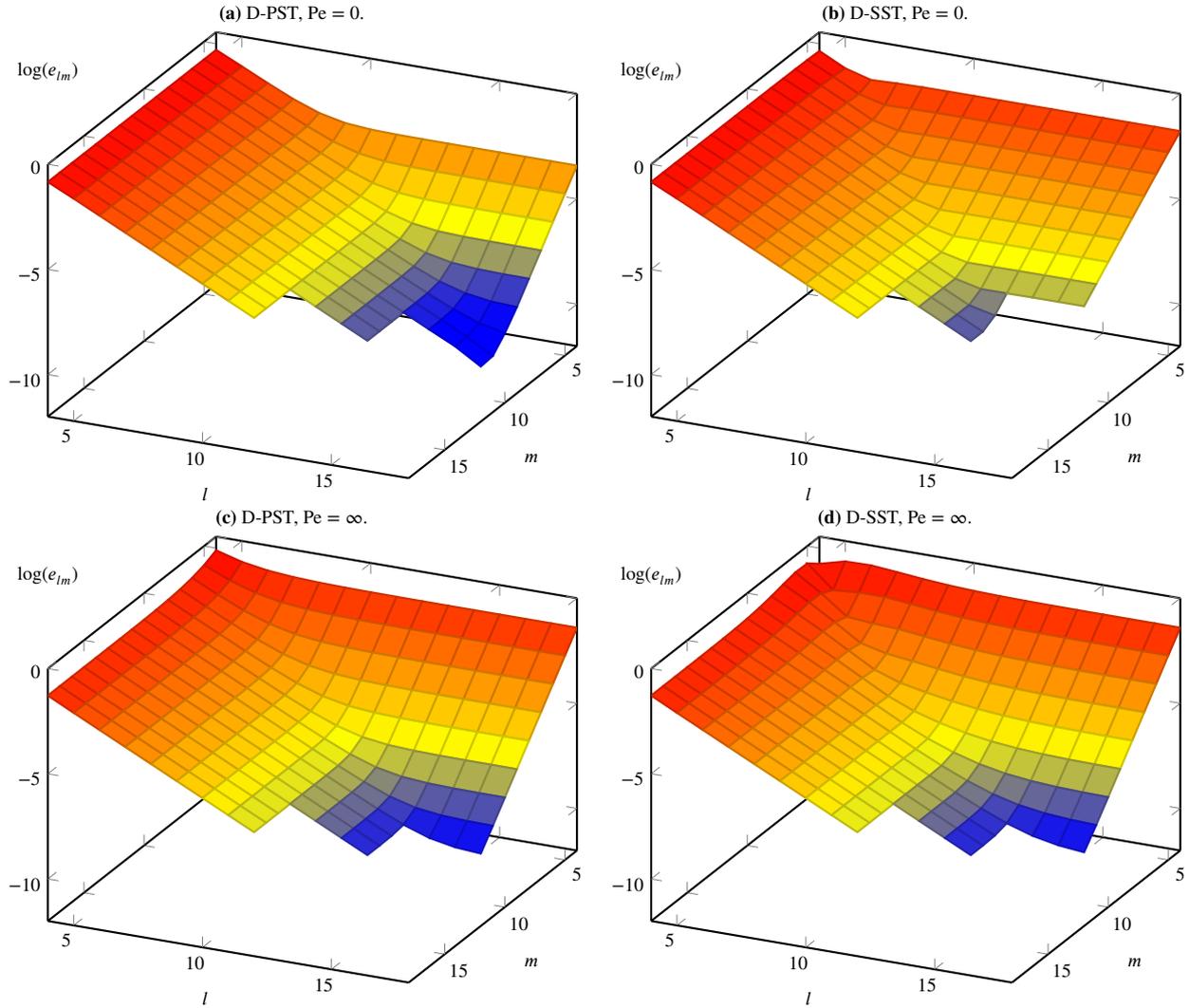

\centering
%\subcaptionbox{FST, advection diffusion case $a = 1$, $\nu = \num{0.01}$ \label{fig:convFSTNuHundred}}{\logPlotThreeD{FSTQuadNu0.01}}%
%\subcaptionbox{SST, advection diffusion case $a = 1$, $\nu = \num{0.01}$ \label{fig:convSSTNuHundred}}{\logPlotThreeD{SST1DNu0.01}} \\
\subcaptionbox{D-PST,  $\Peclet = 0$. \label{fig:convFSTPeZero}}{\logPlotThreeD{FSTQuadPe0}}%
\subcaptionbox{D-SST,  $\Peclet = 0$. \label{fig:convSSTPeZero}}{\logPlotThreeD{SST1DPe0}} \\
\subcaptionbox{D-PST,  $\Peclet = \infty$. \label{fig:convFSTNuZero}}{\logPlotThreeD{FSTQuadNu0.0}}%
\subcaptionbox{D-SST,  $\Peclet = \infty$. \label{fig:convSSTNuZero}}{\logPlotThreeD{SST1DNu0.0}} \\
\caption[Convergence of $L_2$ error for IBVP 1.]{Convergence visualization of $L_2$ error for model problem IBVP 1.}
\label{fig:convIBVP1}
\end{figure} 

As before, the results of the space-time convergence studies are visualized in convergence surfaces. Four representative convergence surfaces are shown in Figure~\ref{fig:convIBVP1}. The complete set of twelve surfaces can be found in Figure~\ref{fig:appADC1L2} in the Appendix. For all parameter sets, a continuous surface is obtained. Furthermore, for both space-time discretizations, the advection-diffusion cases with increasing P\'{e}clet number present a smooth transition from the pure diffusion to the pure advection case. Comparing the parabolic cases (Figure~\ref{fig:convFSTPeZero} and Figure~\ref{fig:convSSTPeZero}), the convergence surfaces of D-PST and D-SST clearly differ. D-PST reaches smaller error values due to the superconvergence of the discretization with tensor-product elements (Equation~\eqref{eq:thomee}). The hyperbolic cases (Figure~\ref{fig:convFSTNuZero} and Figure~\ref{fig:convSSTNuZero}) show only a slight difference for the simulations on coarse meshes with $\Delta x= \Delta t$, yet, on the finer meshes D-PST and D-SST arrive at very similar results.

\newcommand{\logPlotTwoDSpace}[1]{
%\tikzsetnextfilename{spacetime-Space#1}
\begin{tikzpicture}
\begin{loglogaxis}[
       	height=\plthgtC, width=\pltwdtC,
     	xlabel={$\Delta x$},
	small e axis style,
	cycle list name=myDivPlotCycleList,
        legend style={
               at={(0.5,1.2)},
               anchor=north,
               cells={anchor=west},
               legend columns=4,
       },
       legend entries={$\Peclet = 0$, $\Peclet = 1$, $\Peclet = 10$, $\Peclet = 100$, $\Peclet = 1000$, $\Peclet = \infty$,$\mathcal{O}(\Delta x^2)$},
       legend to name=twoDSpaceLegend#1,
]
	\addplot+[dataline] table [x expr=2/(2^(\coordindex+3)), y index=0] {Data/#1Pe0RealL2ErrorSpace.txt};
	\addplot+[dataline] table [x expr=2/(2^(\coordindex+3)), y index=0] {Data/#1Nu0.1RealL2ErrorSpace.txt};
	\addplot+[dataline] table [x expr=2/(2^(\coordindex+3)), y index=0] {Data/#1Nu0.01RealL2ErrorSpace.txt};
	\addplot+[dataline] table [x expr=2/(2^(\coordindex+3)), y index=0] {Data/#1Nu0.001RealL2ErrorSpace.txt};
	\addplot+[dataline] table [x expr=2/(2^(\coordindex+3)), y index=0] {Data/#1Nu0.0001RealL2ErrorSpace.txt};
	\addplot+[dataline] table [x expr=2/(2^(\coordindex+3)), y index=0] {Data/#1Nu0.0RealL2ErrorSpace.txt};
	\addplot[very thick, dashed, domain=0.001:0.25, samples=100]{(x)^(2)/10};
\end{loglogaxis}
\end{tikzpicture}
}

\begin{figure}
\centering
\parbox{.49\linewidth}{
\flushright \parbox{\pltwdtC}{\centering \subcaption{D-PST, $\Delta t = \frac{1}{65536}$.}}
\hfill \logPlotTwoDSpace{FSTQuad}
}\,\parbox{.49\linewidth}{
\flushright \parbox{\pltwdtC}{\centering \subcaption{D-SST, $\Delta t = \frac{1}{65536}$.}}
\hfill \logPlotTwoDSpace{SST1D}}
%\ref{twoDSpaceLegendFSTQuad}
\includegraphics{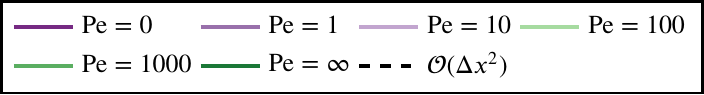}
\caption[Spatial convergence for IBVP 1.]{Spatial convergence for IBVP 1 for six model parameter sets.}
\label{fig:twoDSpace}
\end{figure} 

After this brief view on the convergence surfaces, we now analyze spatial and temporal convergence by means of line plots. Spatial convergence results are presented in Figure~\ref{fig:twoDSpace}. Both space-time methods converge quadratically with respect to $\Delta x$ for the complete model parameter range from $\Peclet = 0$ up to $\Peclet = \infty$ and for all values of $\Delta x$. Besides the constant convergence rates, there is an influence of the P\'{e}clet number on the actual relative error values. The solutions for the more diffusive cases, are slightly more accurate.

\newcommand{\logPlotTwoDTimeFST}[1]{
\tikzsetnextfilename{spacetime-Time#1}
\begin{tikzpicture}
\begin{loglogaxis}[
       height=\plthgtC, width=\pltwdtC,
                  xlabel={$\Delta t$},
		small e axis style,
			cycle list name=myDivPlotCycleList,
]
	\addplot+[dataline] table [x expr=2/(2^(\coordindex+3)), y index=0] {Data//#1Pe0RealL2ErrorTime.txt};
	\addplot+[dataline] table [x expr=2/(2^(\coordindex+3)), y index=0] {Data//#1Nu0.1RealL2ErrorTime.txt};
	\addplot+[dataline] table [x expr=2/(2^(\coordindex+3)), y index=0] {Data//#1Nu0.01RealL2ErrorTime.txt};
	\addplot+[dataline] table [x expr=2/(2^(\coordindex+3)), y index=0] {Data//#1Nu0.001RealL2ErrorTime.txt};
	\addplot+[dataline] table [x expr=2/(2^(\coordindex+3)), y index=0] {Data//#1Nu0.0001RealL2ErrorTime.txt};
	\addplot+[dataline] table [x expr=2/(2^(\coordindex+3)), y index=0] {Data//#1Nu0.0RealL2ErrorTime.txt};
        \addplot[very thick, dotted, domain=0.001:0.25, samples=100]{(x)^(3)};
        \addplot[very thick, dotted, domain=0.005:0.25, samples=100]{(x)^(3)/100};
\end{loglogaxis}
\end{tikzpicture}
}

\newcommand{\logPlotTwoDTimeSST}[1]{
\tikzsetnextfilename{spacetime-Time#1}
\begin{tikzpicture}
\begin{loglogaxis}[
       height=\plthgtC, width=\pltwdtC,
           xlabel={$\Delta t$},
	small e axis style,
		cycle list name=myDivPlotCycleList,
               legend style={
               at={(0.5,1.25)},
               anchor=north,
               cells={anchor=west},
               legend columns=4,
       },
       legend entries={$\Peclet = 0$, $\Peclet = 1$, $\Peclet = 10$, $\Peclet = 100$, $\Peclet = 1000$, $\Peclet = \infty$, $\mathcal{O}(\Delta t^2)$, $\mathcal{O}(\Delta t^3)$},
       legend to name=twoDTimeLegend,
]
	\addplot+[dataline] table [x expr=2/(2^(\coordindex+3)), y index=0] {Data//#1Pe0RealL2ErrorTime.txt};
	\addplot+[dataline] table [x expr=2/(2^(\coordindex+3)), y index=0] {Data//#1Nu0.1RealL2ErrorTime.txt};
	\addplot+[dataline] table [x expr=2/(2^(\coordindex+3)), y index=0] {Data//#1Nu0.01RealL2ErrorTime.txt};
	\addplot+[dataline] table [x expr=2/(2^(\coordindex+3)), y index=0] {Data//#1Nu0.001RealL2ErrorTime.txt};
	\addplot+[dataline] table [x expr=2/(2^(\coordindex+3)), y index=0] {Data//#1Nu0.0001RealL2ErrorTime.txt};
	\addplot+[dataline] table [x expr=2/(2^(\coordindex+3)), y index=0] {Data//#1Nu0.0RealL2ErrorTime.txt};
	\addplot[very thick, dashed, domain=0.001:0.1, samples=100]{(x)^(2)};
        \addplot[very thick, dotted, domain=0.001:0.1, samples=100]{(x)^(3)};
\end{loglogaxis}
\end{tikzpicture}
}

\begin{figure}
\centering
\parbox{.49\linewidth}{
\flushright \parbox{\pltwdtC}{\centering \subcaption{D-PST, $\Delta x = \frac{1}{65536}$.\label{fig:convTwoDFSTTime}}}
\hfill \logPlotTwoDTimeFST{FSTQuad}
}\,\parbox{.49\linewidth}{
\flushright \parbox{\pltwdtC}{\centering \subcaption{D-SST, $\Delta x = \frac{1}{65536}$.\label{fig:convTwoDSSTTime}}}
\hfill \logPlotTwoDTimeSST{SST1D}}
%\ref{twoDTimeLegend}\\
\includegraphics{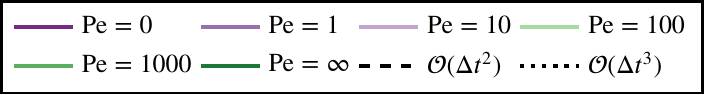}
\caption[Temporal convergence for IBVP 1.]{Temporal convergence for IBVP 1 for six model parameter sets. \label{fig:twoDTime}}
\end{figure} 

Figure~\ref{fig:twoDTime} collects the results of the temporal convergence study. For this model problem, D-PST is observed to converge cubically with respect to the time step for the complete range of P\'{e}clet numbers, despite the fact that the curve of $\Peclet = 0$ is shifted to smaller error values as shown in Figure~\ref{fig:convTwoDFSTTime}. This behavior is in line with the results obtained by Shakib and Hughes in a Fourier analysis of the purely advective and purely diffusive limiting case of this model problem~\cite{shakib1991new}. Note that the specific mesh connectivity (stencil) of the D-PST discretization is used in the Fourier analysis and the results hence do not apply to a D-SST discretization. The D-SST results presented in Figure~\ref{fig:convTwoDSSTTime} show a strong influence of the P\'{e}clet number. The method is second-order time accurate in the parabolic case and third-order accurate in the hyperbolic case. For advection-diffusion cases, we observe a smooth transition of the convergence behavior from second to third order. However, rather than converging at a constant intermediate rate, D-SST converges for the advection-diffusion cases cubically up to some $\Delta t _{\mathrm{turn}} (\Peclet)$, where the convergence rate transitions to two. For smaller P\'{e}clet numbers, the transition occurs at larger time steps, which is earlier in the convergence history.  

\newcommand{\logPlotThreeDNe}[1]{
\tikzsetnextfilename{spacetime-Ne#1}
\begin{tikzpicture}
\begin{axis}
[
height=\plthgtTD, width=\pltwdtTD,
point meta min = -9,
zmax=0,
zmin=-12,
unbounded coords=jump,
y dir=reverse,
xlabel=$l$, ylabel=$m$, 
zlabel=$\log (E_{lm})$,
z label style={at={(axis description cs:0.0,0.95)},,rotate=-90,anchor=north},
]
\addplot3 [surf, mesh/rows=15, mesh/ordering=y varies, z filter/.expression= {z==-15 ? nan : z}] table [col sep=comma]  {Data/#1L2ErrorCoords.csv};
\end{axis}
\end{tikzpicture}
}

\begin{figure}
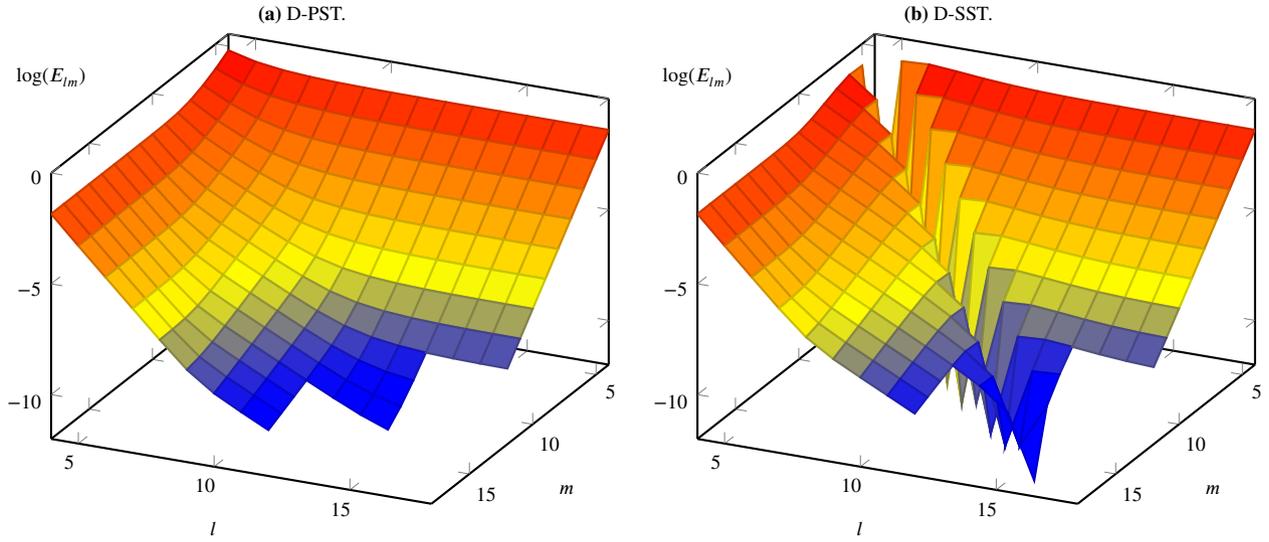

\centering
\subcaptionbox{D-PST. \label{fig:IBVP1FST}}{\logPlotThreeDNe{FSTQuadNu0.0}}%
\subcaptionbox{D-SST. \label{fig:IBVP1SST}}{\logPlotThreeDNe{SST1DNu0.0}} \\
\caption[Convergence of nodal error for IBVP 1 with $\Peclet=\infty$.]{Convergence visualization of nodal error measure for IBVP 1 with $\Peclet=\infty$.}
\label{fig:convIBVP1Ne}
\end{figure} 

So far, we discussed the convergence of the $L_2$ errors, but, also the convergence surfaces of the nodal error measure (Equation~\eqref{eq:adNe}) show interesting features of the discretization methods. The complete set of the twelve nodal error based convergence surfaces can be found in Figure~\ref{fig:appADC1Ne} in the Appendix~\ref{sec:appendix}. The results of the hyperbolic case are presented in Figure~\ref{fig:convIBVP1Ne}. The nodal error visualization of the D-PST results (Figure~\ref{fig:IBVP1FST}) is a continuous surface as for the $L_2$ error. In contrast, the D-SST results (Figure~\ref{fig:IBVP1SST}) show a strong discontinuity for the simulations with $\Delta t = \Delta x$. For these cases, the space-time finite element edges align with the characteristic curves along which the solution is transported. We observe that the finite element approximation coincides with the exact solution at the nodes (up to a round off error $\epsilon < \num{1.0e-10}$) for all refinement levels. This astonishing behavior is described by Demkowicz and Oden as `extra superconvergence'~\cite{demkowicz1986adaptive}. Away from the diagonal $l=m$, error values are obtained that are similar to the ones of D-PST.

\newcommand{\plthgtCNew}{0.2*\textheight}
\newcommand{\pltwdtCNew}{0.37*\textwidth}
\newcommand{\logPlotTwoDSpaceFSTNe}[1]{
\tikzsetnextfilename{spacetime-SpaceNe#1}
\begin{tikzpicture}
\begin{loglogaxis}[
       height=\plthgtCNew, width=\pltwdtC,
       scale only axis,
    	xlabel={$\Delta x, \left( \Delta t = 1/65536\right)$},
       	ymax = 0.5,
	ymin = 1e-14,
       	ylabel= {$E_{lm}$},
	y label style={at={(axis description cs:0.0,1.0)},,rotate=-90,anchor=south east},
	grid=major,
	cycle list name=myDivPlotCycleList,
	legend style={
               at={(0.5,1.2)},
               anchor=north,
               cells={anchor=west},
               legend columns=4,
       },
       legend entries={$\Peclet = 0$, $\Peclet = 1$, $\Peclet = 10$, $\Peclet = 100$, $\Peclet = 1000$, $\Peclet = \infty$, $\mathcal{O}(\Delta t^2) \,\text{or}\, \mathcal{O}(\Delta x^2)$, $\mathcal{O}(\Delta x^4)$},
       legend to name=twoDSpaceTimeNeLegend,
]
	\addplot+[dataline] table [x expr=2/(2^(\coordindex+3)), y index=0] {Data/#1Pe0L2ErrorSpace.txt};
	\addplot+[dataline] table [x expr=2/(2^(\coordindex+3)), y index=0] {Data/#1Nu0.1L2ErrorSpace.txt};
	\addplot+[dataline] table [x expr=2/(2^(\coordindex+3)), y index=0] {Data/#1Nu0.01L2ErrorSpace.txt};
	\addplot+[dataline] table [x expr=2/(2^(\coordindex+3)), y index=0] {Data/#1Nu0.001L2ErrorSpace.txt};
	\addplot+[dataline] table [x expr=2/(2^(\coordindex+3)), y index=0] {Data/#1Nu0.0001L2ErrorSpace.txt};
	\addplot+[dataline] table [x expr=2/(2^(\coordindex+3)), y index=0] {Data/#1Nu0.0L2ErrorSpace.txt};
        \addplot[very thick, dashed, domain=0.001:0.25, samples=100]{(x)^(2)/3};
        \addplot[very thick, dashdotted, domain=0.001:0.25, samples=100]{(x)^(4)};
\end{loglogaxis}
\end{tikzpicture}
}

\newcommand{\logPlotTwoDSpaceTimeSSTNe}[1]{
\tikzsetnextfilename{spacetime-SpaceTimeNe#1}
\begin{tikzpicture}
\begin{loglogaxis}[
        height=\plthgtCNew, width=\pltwdtC,
        scale only axis,
    	xlabel={$\Delta t = \Delta x$},
       	ymax = 0.5,
	ymin = 1e-14,
       	ylabel= {$E_{lm}$},
	y label style={at={(axis description cs:0.0,1.0)},,rotate=-90,anchor=south east},
	grid=major,
	cycle list name=myDivPlotCycleList,
	]
	\addplot+[dataline] table [x expr=2/(2^(\coordindex+3)), y index=0] {Data/#1Pe0L2ErrorSpaceTime.txt};
	\addplot+[dataline] table [x expr=2/(2^(\coordindex+3)), y index=0] {Data/#1Nu0.1L2ErrorSpaceTime.txt};
	\addplot+[dataline] table [x expr=2/(2^(\coordindex+3)), y index=0] {Data/#1Nu0.01L2ErrorSpaceTime.txt};
	\addplot+[dataline] table [x expr=2/(2^(\coordindex+3)), y index=0] {Data/#1Nu0.001L2ErrorSpaceTime.txt};
	\addplot+[dataline] table [x expr=2/(2^(\coordindex+3)), y index=0] {Data/#1Nu0.0001L2ErrorSpaceTime.txt};
	\addplot+[dataline] table [x expr=2/(2^(\coordindex+3)), y index=0] {Data/#1Nu0.0L2ErrorSpaceTime.txt};
	\addplot[very thick, dashed,domain=0.000125:0.25, samples=100]{(x)^(2)/3};
\end{loglogaxis}
\end{tikzpicture}
}

\begin{figure}
\centering
\parbox{.49\linewidth}{\flushright \parbox{\pltwdtC}{\centering \subcaption{Spatial convergence of D-PST. \label{fig:convTwoDFSTSpaceNe}}}
 \logPlotTwoDSpaceFSTNe{FSTQuad}}\,\parbox{.49\linewidth}{
\flushright \parbox{\pltwdtC}{\centering \subcaption{Space-time convergence of D-SST. \label{fig:convTwoDSSTSpaceTimeNe}}}
\logPlotTwoDSpaceTimeSSTNe{SST1D}}
%\ref{twoDSpaceTimeNeLegend}\\
\includegraphics{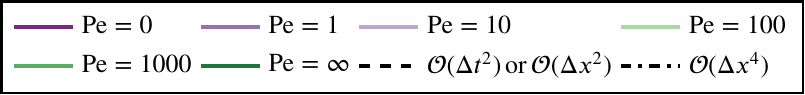}
\caption[Convergence of nodal error for IBVP 1 for six parameter sets.]{Convergence in nodal error measure for IBVP 1 for six parameter sets.\label{fig:twoDSpaceTimeNe}}
\end{figure} 

Spatial convergence results of the D-PST method in nodal error measure are extracted as line plots and shown in Figure~\ref{fig:convTwoDFSTSpaceNe}. We see once more a strong influence of the P\'{e}clet number on the convergence behavior. In the pure advection case, the method converges in the nodal error measure with fourth order up to a $\Delta x_{\mathrm{turn}} (\Peclet)$ and then transitions to second order. For smaller P\'{e}clet numbers, the transition occurs at larger $\Delta x$, i. e., earlier in the convergence history. In the simulations with the fourth order convergence relative to each other, the element size in time direction $\Delta t$ is very small compared to $\Delta x$. In consequence, the small $\Delta t$ leads to such a small stabilization parameter $\tau_{\mathrm{SUPG}}$ (Equation~\eqref{eq:tauAD}), that the influence of the SUPG term vanishes and the Galerkin method is recovered. The nodal error of the Galerkin method for the pure advection case is fourth order accurate with respect to $\Delta x$, as shown in the Fourier analysis of Shakib and Hughes~\cite{shakib1991new}. In the pure diffusion case, the method is second order accurate with respect to $\Delta x$ over the entire element size range.

Figure~\ref{fig:convTwoDSSTSpaceTimeNe} shows nodal error results of the D-SST method for the six model parameter sets along the space-time diagonal $\Delta t = \Delta x$. Most notable is the `extra superconvergence` of the SST method with characteristics aligned element edges for the pure advection case. For the five other model parameter sets, we observe a second order space-time convergence along the diagonal $\Delta t = \Delta x$. The curves of the cases with $\Peclet = 0 $ and $\Peclet = 1$ essentially coincide, while the other cases show smaller error values for higher P\'{e}clet numbers.

Unfortunately, it is highly unlikely that the finite element edges of higher-dimensional space-time meshes are aligned with the solution characteristics for general flow conditions. Therefore, we come to the following outlook for problems of engineering interest. Under the assumption that our findings carry over from the scalar one-dimensional advection-diffusion cases to higher-dimensional cases modeled with (in)compressible Navier--Stokes equations, we expect a spatial accuracy of second order and a temporal accuracy between second and third order. As we have seen, the temporal accuracy of the time-discontinuous space-time methods tends towards third-order depending on how advection-dominated the test case is, on the element type used for discretization, and on the time-(in)dependence and treatment of the boundary conditions.

\section{Piston Ring Test Case}
\label{sec:prp}
The time-discontinuous space-time discretizations (D-PST and D-SST) have several advantages, e.g., with tensor-product elements superconvergence of the $L_2$-error at the final time can be achieved. However, in practice superconvergence is hard to obtain as it is contingent on several factors, e.g., the treatment of time-dependent boundary conditions. Therefore, we use in this section the more flexible time-continuous simplex space-time discretization (C-SST). 

\begin{figure}
\centering
\includegraphics[width=10cm,trim={0cm 0cm 0cm 0cm},clip]{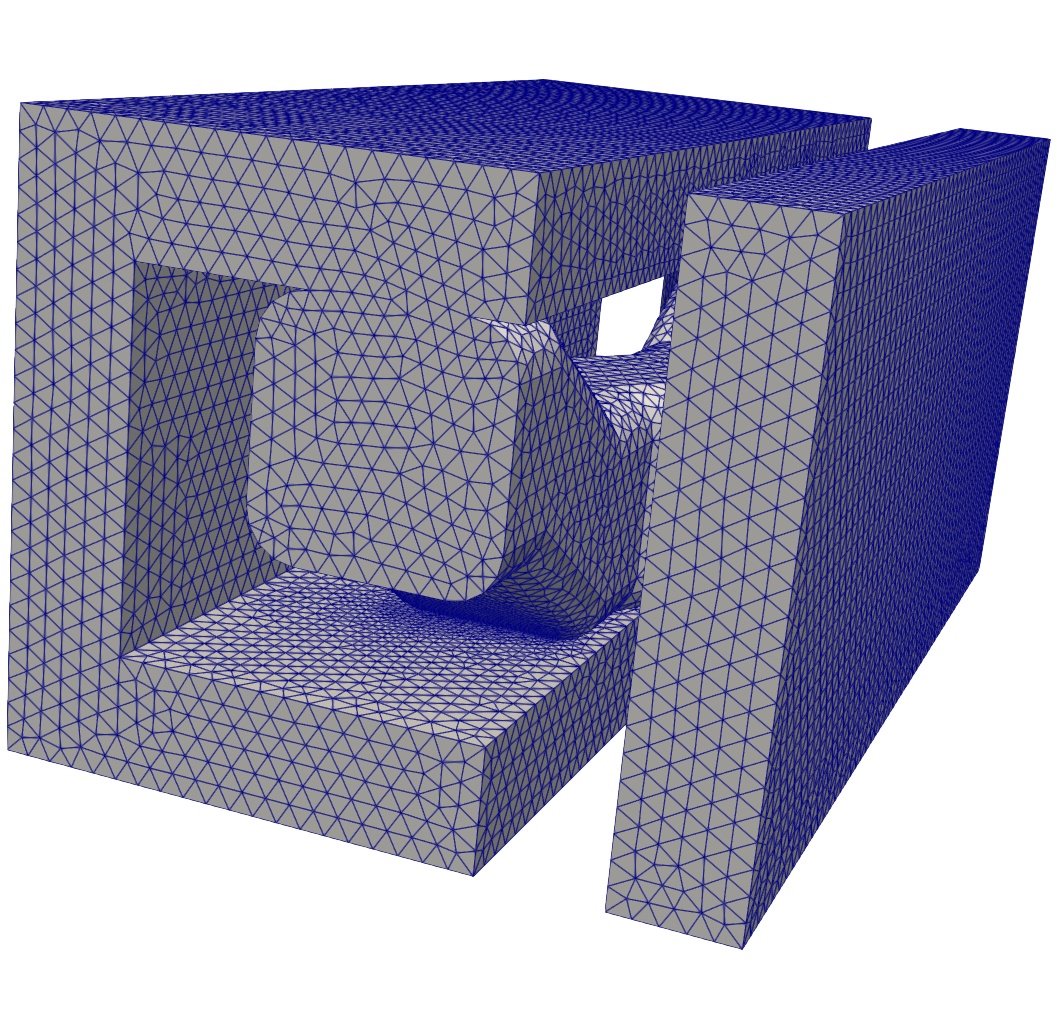}
\caption{Tetrahedral space-time mesh.}
\label{fig:prpMesh}
\end{figure}

The purpose of the following test case is to demonstrate the capability of time-continuous simplex space-time discretizations (C-SST) to account for complex changes of the spatial computational domain. In this particular simulation, the connectivity of the spatial domain changes multiple times. The boundary conforming tetrahedral space-time mesh is shown in Figure~\ref{fig:prpMesh}.

As a motivational example, we consider the \prp on an internal combustion engine. The piston rings are employed to seal the high-pressure gas in the combustion chamber (i), to prevent engine oil from leaking into the combustion chamber (ii), and to dissipate heat from the piston to the surrounding cylinder to prevent overheating of the piston (iii). In this test case, we investigate the heat flux in a simplified model of a \prp. 

\renewcommand{\swW}{0.24\textheight}
\newcommand{\myY}{1.25}
\newcommand{\myL}{0.25}
\newcommand{\myX}{0.6}

\begin{figure}
\centering 
\begin{tikzpicture}[scale=5,
    axis/.style={thick, ->},
    every node/.style={color=black}
    ]
       % Include image
       \node at (-5*\myL,-\myL) [inner sep=0pt, anchor=south west] {\includegraphics[width=7.5cm,trim={0cm 0cm 0cm 0cm},clip]{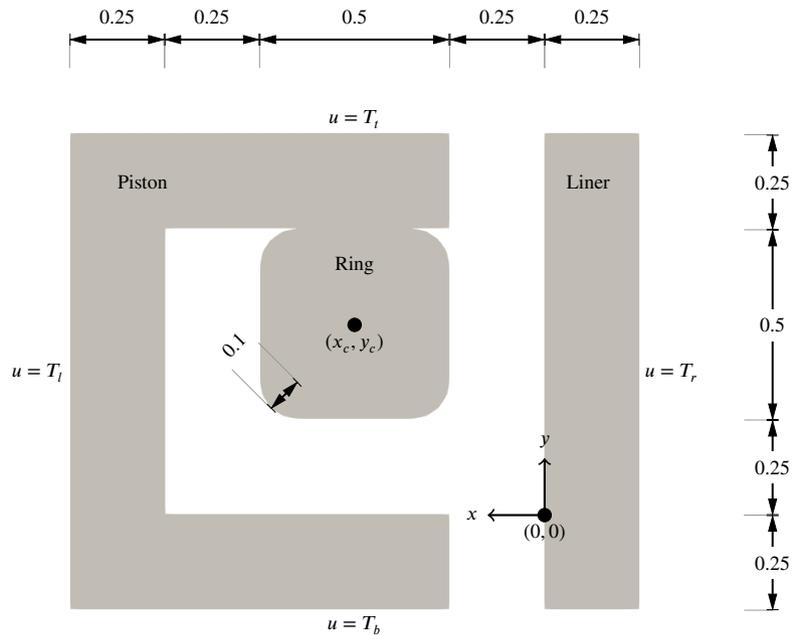}};
       
       % Label parts
       \node at (-4.5*\myL,3.5*\myL) [inner sep=0pt, anchor=west] {Piston};
       \node at (-2*\myL,2.75*\myL) [inner sep=1pt, anchor=north] {Ring};
       \node at (0.225*\myL,3.5*\myL) [inner sep=0pt, anchor=west] {Liner};
       
       % Coordinates
       \draw[fill] (-2*\myL,2*\myL) circle (0.5pt);
       \node at (-2*\myL,2*\myL) [inner sep=3pt, anchor=north] {$(x_c,y_c)$}; 
       \draw[fill] (0,0) circle (0.5pt);
       \node at (0,0) [inner sep=3pt, anchor=north] {$(0,0)$};
       
       % Cosy
        \draw[axis] (0,0) -- (-0.15,0) node(xline) [left] {$x$};
       \draw[axis] (0,0) -- (0,0.15) node(yline) [above] {$y$};
       
       % Temperatures
       \node at (-2*\myL,4*\myL) [inner sep=2pt, anchor=south] {$u = T_t$};
        \node at (-2*\myL,-\myL) [inner sep=2pt, anchor=north] {$u = T_b$};
        \node at (-5*\myL,1.5*\myL) [inner sep=2pt, anchor=east] {$u = T_l$};
        \node at (\myL,1.5*\myL) [inner sep=2pt, anchor=west] {$u = T_r$};
       
       % length in vertical direction
       \dimline[line style={line width=0.75pt}, label style={rotate=-90},extension start length=-0.3,extension end length=-0.3]{(\myX,-\myL)}{(\myX,0)}{$0.25$};
       \dimline[line style={line width=0.75pt},label style={rotate=-90},extension start length=-0.3,extension end length=-0.3]{(\myX,0)}{(\myX,\myL)}{$0.25$};
       \dimline[line style={line width=0.75pt},label style={rotate=-90},extension start length=0,extension end length=0]{(\myX,\myL)}{(\myX,3*\myL)}{$0.5$};
       \dimline[line style={line width=0.75pt},label style={rotate=-90},extension start length=-0.3,extension end length=-0.3]{(\myX,3*\myL)}{(\myX,4*\myL)}{$0.25$};
       
       % length in horizontal direction
       \dimline[line style={line width=0.75pt},label style={above=0.5ex}, extension start length=0.3,extension end length=0.3]{(-5*\myL,\myY)}{(-4*\myL,\myY)}{$0.25$};
       \dimline[line style={line width=0.75pt},label style={above=0.5ex}, extension start length=0.3,extension end length=0.3]{(-4*\myL,\myY)}{(-3*\myL,\myY)}{$0.25$};
       \dimline[line style={line width=0.75pt},label style={above=0.5ex}, extension start length=0,extension end length=0]{(-3*\myL,\myY)}{(-\myL,\myY)}{$0.5$};
       \dimline[line style={line width=0.75pt},label style={above=0.5ex}, extension start length=0.3,extension end length=0.3]{(-\myL,\myY)}{(0,\myY)}{$0.25$};
       \dimline[line style={line width=0.75pt},label style={above=0.5ex}, extension start length=0.3,extension end length=0.3]{(0,\myY)}{(\myL,\myY)}{$0.25$};
       
       % Radius ring corners
       \dimline[line style={line width=0.75pt},label style={above=4.5ex}, extension start length=2,extension end length=2]{(-0.65,0.35)}{(-0.72,0.28)}{$0.1$};
       
\end{tikzpicture}
\caption[Setup of piston ring test case.]{Setup of piston ring test case.}
\label{fig:prpSetup}
\end{figure}

Figure~\ref{fig:prpSetup} shows the two-dimensional geometry of a schematic \prp with only one ring. The considered geometry includes a part of the piston around the groove in which the piston ring is located, as well as a part of the cylinder liner which comes into contact with the piston ring. As shown in Figure~\ref{fig:prpSetup}, the piston ring is represented by a square with a generic side length of 0.5; its corners are rounded with radius $r=0.1$.

In the following simulation, we investigate the conductive heat transfer in the metal parts and across the contact interfaces between piston, ring, and liner. The heat transfer is modeled with the parabolic case of Equation~\eqref{eq:ADstrong}, i.e., the advection velocity is set to zero and we obtain a P\'{e}clet number of zero. The thermal diffusivity in the solids $k = \frac{\kappa}{\rho c_p}$, that accounts for the thermal conductivity $\kappa$, density $\rho$, and specific heat $c_p$, is here modeled with a generic diffusion coefficient $k=0.495$, as employed in Equation~\eqref{eq:ADstrong}. The test case is further characterized by the temperatures
\begin{equation} T_t = 423.15,  \quad T_b = 403.15, \quad T_l = 403.15 + 20 \cdot \frac{y+0.25}{1.25}, \quad \text{and} \; T_r = 373.15 \end{equation} prescribed as Dirichlet boundary conditions on the edges as indicated in Figure~\ref{fig:prpSetup}. On all remaining boundaries, homogeneous Neumann boundary conditions are assumed. The initial boundary value problem is completed with the initial condition
\begin{equation} u(x,y,t=0) = T_0 = 373.15 + h(x-0.01) \cdot \left(30+20 \cdot \frac{y+0.25}{1.25}\right), \end{equation} where $h(x)$ denotes the Heaviside function. 

\begin{figure}
\centering
\begin{tikzpicture}
\centering
\begin{axis}[
   height=0.3\textheight, width=0.7\textwidth,
    ylabel={Ring center coordinates},
    ylabel near ticks,
    ytick= {0.25,0.3,0.4,0.5},
   	xlabel= {Time $t$},
		xtick={0,0.2,0.6,0.8,1.2,1.4,1.8,2.0,2.4,2.6},
   	grid=major,
   	major grid style={dotted},
		legend style={
		at={(0.5,0.9)},
		anchor=north,
		cells={anchor=west},
		legend columns=1},
   	legend entries={$x_c(t)$, $y_c(t)$},
]    
	\addplot+ [no markers,black] coordinates {
		(0,   0.5)
		(0.2,0.5)
		(0.6,0.5)
		(0.8,0.5)
		(1.2,0.25)
		(1.4,0.25)
		(1.8,0.25)
		(2.0,0.25)
		(2.4,0.5)
		(2.6,0.5)
		};    
		
			\addplot+ [no markers,black,dashed] coordinates {
		(0,   0.5)
		(0.2,0.5)
		(0.6,0.25)
		(0.8,0.25)
		(1.2,0.25)
		(1.4,0.25)
		(1.8,0.5)
		(2.0,0.5)
		(2.4,0.5)
		(2.6,0.5)
		};    
\end{axis}
\end{tikzpicture}
\caption[Ring motion.]{Ring motion.}
\label{fig:prpRingMotion}
\end{figure}
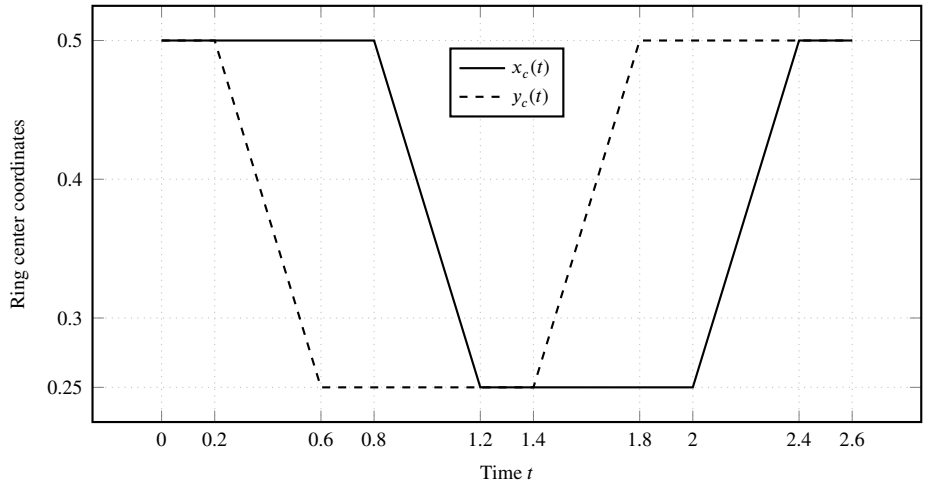

\renewcommand{\swW}{0.24\textheight}

\newcommand{\addTime}[2]{
%\tikzsetnextfilename{spacetime-sw#1}
\begin{tikzpicture}[
    axis/.style={thick, ->},
    every node/.style={color=black}
    ]
       \node at (0,0) [inner sep=0pt, anchor=south west] {\includegraphics[height=\swW,trim={6cm 4cm 6cm 4cm},clip]{#1}};
       \node at (0.135*\swW,0.135*\swW) [inner sep=0pt, anchor=south west]{$t = #2$};
\end{tikzpicture}
}

\begin{figure}
\centering
\vspace{-2ex}
\addTime{Pics/RingAt0.png}{0}\hspace{2em}
\addTime{Pics/RingAt26.png}{2.6}\\\vspace{0.5ex}
\addTime{Pics/RingAt04.png}{0.4}\hspace{2em}
\addTime{Pics/RingAt24.png}{2.4}\\\vspace{0.5ex}
\addTime{Pics/RingAt08.png}{0.8}\hspace{2em}
\addTime{Pics/RingAt20.png}{2.0}\\\vspace{0.5ex}
\addTime{Pics/RingAt12.png}{1.2}\hspace{2em}
\addTime{Pics/RingAt16.png}{1.6}\\
\includegraphics[width=0.3\textwidth]{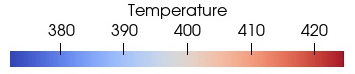}
\caption[Temporal evolution of temperature in \prp.]{Temporal evolution of temperature in \prp.}
\label{fig:temperature}
\end{figure}

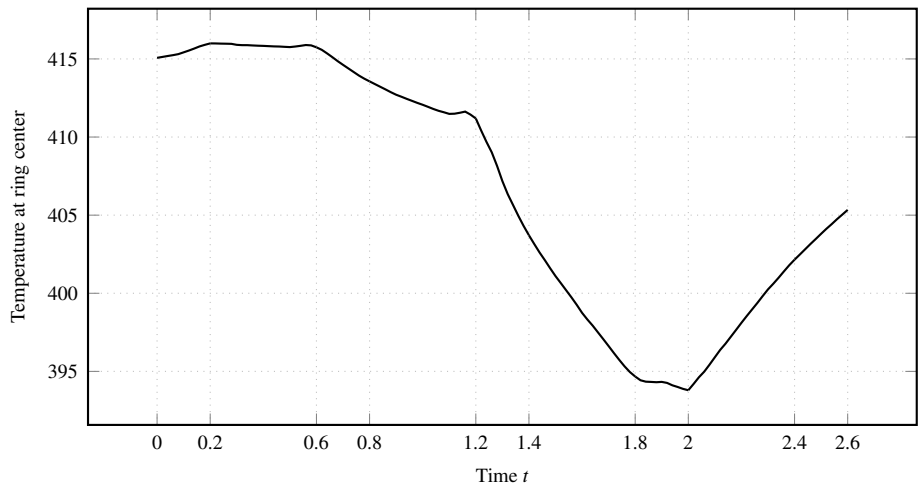
\begin{figure}
\centering
\begin{tikzpicture}
\centering
\begin{axis}[
   height=0.3\textheight, width=0.7\textwidth,
    ylabel={Temperature at ring center},
    ylabel near ticks,
   	xlabel= {Time $t$},
	xtick={0,0.2,0.6,0.8,1.2,1.4,1.8,2.0,2.4,2.6},
   	grid=major,
   	major grid style={dotted}
]    
\addplot [no markers,black] table [x index = 5, y index = 0, col sep=comma, header=true]  {Pics/PRP/temperatureOverTime1.csv};
\addplot [no markers,black] table [x index = 5, y index = 0, col sep=comma, header=true]  {Pics/PRP/temperatureOverTime2.csv};
\addplot [no markers,black] table [x index = 5, y index = 0, col sep=comma, header=true]  {Pics/PRP/temperatureOverTime3.csv};
\addplot [no markers,black] table [x index = 5, y index = 0, col sep=comma, header=true]  {Pics/PRP/temperatureOverTime4.csv};
\addplot [no markers,black] table [x index = 5, y index = 0, col sep=comma, header=true]  {Pics/PRP/temperatureOverTime5.csv};
\addplot [no markers,black] table [x index = 5, y index = 0, col sep=comma, header=true]  {Pics/PRP/temperatureOverTime6.csv};
\addplot [no markers,black] table [x index = 5, y index = 0, col sep=comma, header=true]  {Pics/PRP/temperatureOverTime7.csv};
\addplot [no markers,black] table [x index = 5, y index = 0, col sep=comma, header=true]  {Pics/PRP/temperatureOverTime8.csv};
\addplot [no markers,black] table [x index = 5, y index = 0, col sep=comma, header=true]  {Pics/PRP/temperatureOverTime9.csv};
\end{axis}
\end{tikzpicture}
\caption{Temperature evolution at ring center.}
\label{fig:prpTemperatureRingCenter}
\end{figure}

What makes this test case challenging is the ring motion. In the course of an engine working cycle, the piston ring is in contact with different parts of the piston and the liner. We consider a prescribed ring motion defined by the ring center position, $(x_c(t),y_c(t))$, as shown in Figure~\ref{fig:prpRingMotion}. During the simulated time interval $t \in (0,2.6)$, the ring is first in contact with the upper edge of the piston groove, then moves downwards and is free-floating for $t \in (0.2,0.6)$, before it comes into contact with the lower edge of the piston groove. These three states are also visualized in the first three figures in the left column of Figure~\ref{fig:temperature}. In the following, the ring moves towards the liner, slides upwards along the liner and finally returns to the initial position.

In a C-SST approach~\cite{danwitz2021four}, the given ring motion is included in the computational space-time domain as shown in Figure~\ref{fig:prpMesh}. We used GMSH~\cite{geuzaine2009gmsh} to discretize the domain with a fully unstructured space-time mesh. The resulting mesh consists of 151,911 tetrahedral elements connecting 35,341 nodes. Moreover, the mesh is refined in areas where large spatial and temporal solution gradient are expected, i.e., the curves where the ring comes into contact with the piston and the liner in the course of the simulation. 

The simulation results are collected in Figure~\ref{fig:temperature} and Figure~\ref{fig:prpTemperatureRingCenter}. Figure~\ref{fig:temperature} shows the temperature distribution in the \prp at eight time instances. Most of the time, the temperature solution in the piston and liner parts closely follows the prescribed boundary conditions. Larger spatial temperature variations are primarily encountered in ring. In particular at $t=2.0$, the ring directly connects the hot upper groove edge of the piston with the cooler liner.  As indicated by the large temperature gradients, this configuration leads to the maximal conductive heat transfer.  

In Figure~\ref{fig:prpTemperatureRingCenter}, the temperature at the piston ring center, $T(x_c(t), y_c(t))$, is plotted over time. It is observed that the temperature is approximately constant during the interval $t \in (0.2,0.6)$, which is expected as there is no conductive heat transfer to or from the free-floating ring. The strongest decrease in temperature is observed during the interval $t \in (1.2,1.4)$, when the ring first comes into contact with the cooler liner. The minimal temperature value is reached at $t=2.0$, before the ring again detaches from the liner and is heated from the upper edge of the piston groove.

In summary, the obtained results confirm that C-SST discretizations can easily handle spatial computational domains undergoing complex changes.

\section{Conclusions}
\label{sec:conclusion}
In this paper, we described four space-time finite element methods that result from the combination of tensor-product and simplex-type elements with globally continuous interpolations of the spatial domain and a continuous or discontinuous interpolation in temporal direction. Descriptive naming was proposed, and all four methods were successfully applied to an advection-diffusion model problem. Theoretical background and a detailed numerical convergence analysis were presented for the time-discontinuous space-time methods (D-PST and D-SST). Based on the $L_2$-error at the final time, it was observed that the temporal accuracy of the methods tends towards third-order. For a parabolic model problem, the influence of time-dependent boundary conditions, their treatment, and the element type (prismatic or simplex) of the discretization was studied. For a second model problem with analytical solution, the influence of the element type and the P\'{e}clet number on the convergence behavior was precisely characterized. Moreover, we used the flexible time-continuous simplex space-time (C-SST) method in a challenging heat transfer simulation based on a piston-ring geometry.

\section{Acknowledgment}
\label{sec:acknowledgment}
%This work was supported by the German Federal Ministry for Economic Affairs and Energy through the Central Innovation Programme for Small and Medium Enterprises in the cooperation project ''Simulationstechnik Kolben-Kolbenring-Zylinder'' under grant number KF3462301PO4. 
The authors gratefully acknowledge the computing time granted by the JARA Vergabegremium and provided on the JARA Partition part of the supercomputer CLAIX at RWTH Aachen University.
% Furthermore, we want to thank Philipp Knechtges for his helpful remarks during the preparation of the manuscript.

\clearpage

\bibliography{main}

\begin{thebibliography}{10}

\bibitem{Hirt1974}
Hirt CW, Amsden AA, Cook JL. An arbitrary {L}agrangian-{E}ulerian computing
  method for all flow speeds.  {\it J {C}omput {P}hys. }1974;14(3):227--253.

\bibitem{Liska2010}
Liska R, Shashkov M, Vachal P, Wendroff B. Optimization-based synchronized
  flux-corrected conservative interpolation (remapping) of mass and momentum
  for arbitrary {L}agrangian-{E}ulerian methods.  {\it J {C}omput {P}hys.
  }2010;229(5):1467--1497.

\bibitem{Taylor1937}
Taylor GI, Green AE. Mechanism of the production of small eddies from large
  ones.  {\it P {R}oy {S}oc {L}ond {A} {M}at. }1937;158(895):499--521.
\newblock \url{https://doi.org/10.1098/rspa.1937.0036},
  \url{http://rspa.royalsocietypublishing.org/content/158/895/499}.

\bibitem{Knupp1999}
Knupp PM. Winslow smoothing on two-dimensional unstructured meshes.  {\it Eng
  {C}omput. }1999;15:263--268.

\bibitem{Kamm2000}
Kamm J. {\it Evaluation of the {S}edov-von {N}eumann-{T}aylor blast wave
  solution. } Technical {R}eport LA-UR-00-6055: Los {A}lamos {N}ational
  {L}aboratory; 2000.

\bibitem{Kucharik2003}
Kucharik M, Shashkov M, Wendroff B. An efficient linearity-and-bound-preserving
  remapping method.  {\it J {C}omput {P}hys. }2003;188(2):462--471.

\bibitem{Blanchard2015}
Blanchard G, Loubere R. {\it High-Order {C}onservative {R}emapping with a
  posteriori {MOOD} stabilization on polygonal meshes. }
  \url{https://hal.archives-ouvertes.fr/hal-01207156}, the {HAL} {O}pen
  {A}rchive, hal-01207156. Accessed January 13, 2016; 2015.

\bibitem{Burton2013}
Burton DE, Kenamond MA, Morgan NR, Carney TC, Shashkov MJ. An intersection
  based {ALE} scheme {(xALE)} for cell centered hydrodynamics {(CCH)}.  In:
  Talk at {M}ultimat 2013, {I}nternational {C}onference on {N}umerical
  {M}ethods for {M}ulti-{M}aterial {F}luid {F}lows; September 2--6, 2013; San
  {F}rancisco.
\newblock LA-UR-13-26756.2.

\bibitem{Berndt2011}
Berndt M, Breil J, Galera S, Kucharik M, Maire PH, Shashkov M. Two-step hybrid
  conservative remapping for multimaterial arbitrary {L}agrangian-{E}ulerian
  methods.  {\it J {C}omput {P}hys. }2011;230(17):6664--6687.

\bibitem{Kucharik2012}
Kucharik M, Shashkov M. One-step hybrid remapping algorithm for multi-material
  arbitrary {L}agrangian-{E}ulerian methods.  {\it J {C}omput {P}hys.
  }2012;231(7):2851--2864.

\bibitem{Breil2015}
Breil J, Alcin H, Maire PH. A swept intersection-based remapping method for
  axisymmetric {ReALE} computation.  {\it Int {J} {N}umer {M}eth {F}l.
  }2015;77(11):694--706.
\newblock Fld.3996.

\bibitem{Barth1997}
Barth TJ. Numerical methods for gasdynamic systems on unstructured meshes.  In:
   Kroner D, Rohde C, Ohlberger M, eds. {\it An {I}ntroduction to {R}ecent
  {D}evelopments in {T}heory and {N}umerics for {C}onservation {L}aws,
  {P}roceedings of the {I}nternational {S}chool on {T}heory and {N}umerics for
  {C}onservation {L}aws}, Lecture {N}otes in {C}omputational {S}cience and
  {E}ngineering. Berlin: Springer 1997.
\newblock ISBN 3-540-65081-4.

\bibitem{Lauritzen2011}
Lauritzen P, Erath C, Mittal R. On simplifying `incremental remap'-based
  transport schemes.  {\it J {C}omput {P}hys. }2011;230(22):7957--7963.

\bibitem{Klima2017}
Klima M, Kucharik M, Shashkov M. Local error analysis and comparison of the
  swept- and intersection-based remapping methods.  {\it Commun {C}omput
  {P}hys. }2017;21(2):526--558.

\bibitem{Dukowicz2000}
Dukowicz JK, Baumgardner JR. Incremental remapping as a transport/advection
  algorithm.  {\it J {C}omput {P}hys. }2000;160(1):318--335.

\bibitem{Kucharik2011}
Kucharik M, Shashkov M. Flux-based approach for conservative remap of
  multi-material quantities in {2D} arbitrary {L}agrangian-{E}ulerian
  simulations.  In:  Fo\v{r}t J, F{\"{u}}rst J, Halama J, Herbin R, Hubert F,
  eds. {\it Finite {V}olumes for {C}omplex {A}pplications {VI} {P}roblems \&
  {P}erspectives},  Springer {P}roceedings in {M}athematics, vol. 1: Springer
  2011 (pp. 623--631).

\bibitem{Kucharik2014}
Kucharik M, Shashkov M. Conservative multi-material remap for staggered
  multi-material arbitrary {L}agrangian-{E}ulerian methods.  {\it J {C}omput
  {P}hys. }2014;258:268--304.

\bibitem{Loubere2005}
Loubere R, Shashkov M. A subcell remapping method on staggered polygonal grids
  for arbitrary-{L}agrangian-{E}ulerian methods.  {\it J {C}omput {P}hys.
  }2005;209(1):105--138.

\bibitem{Caramana1998}
Caramana EJ, Shashkov MJ. Elimination of artificial grid distortion and
  hourglass-type motions by means of {L}agrangian subzonal masses and
  pressures.  {\it J {C}omput {P}hys. }1998;142(2):521--561.

\bibitem{Hoch2009}
Hoch P. {\it An arbitrary {L}agrangian-{E}ulerian strategy to solve
  compressible fluid flows. } Technical {R}eport: CEA; 2009.
\newblock HAL: hal-00366858.
  https://hal.archives-ouvertes.fr/docs/00/36/68/58/PDF/ale2d.pdf. Accessed
  January 13, 2016.

\bibitem{Shashkov1996}
Shashkov M. {\it Conservative {F}inite-{D}ifference {M}ethods on {G}eneral
  {G}rids}.
\newblock Boca Raton, Florida: CRC {P}ress; 1996.
\newblock ISBN 0-8493-7375-1.

\bibitem{Benson1992}
Benson DJ. Computational methods in {L}agrangian and {E}ulerian hydrocodes.
  {\it Comput {M}ethod {A}ppl {M}. }1992;99(2--3):235--394.

\bibitem{Margolin2003}
Margolin LG, Shashkov M. Second-order sign-preserving conservative
  interpolation (remapping) on general grids.  {\it J {C}omput {P}hys.
  }2003;184(1):266--298.

\bibitem{Kenamond2013}
Kenamond MA, Burton DE. Exact intersection remapping of multi-material
  domain-decomposed polygonal meshes.  In: Talk at {M}ultimat 2013,
  {I}nternational {C}onference on {N}umerical {M}ethods for {M}ulti-{M}aterial
  {F}luid {F}lows; September 2--6, 2013; San {F}rancisco.
\newblock LA-UR-13-26794.

\bibitem{Dukowicz1984}
Dukowicz J. Conservative rezoning (remapping) for general quadrilateral meshes.
   {\it J {C}omput {P}hys. }1984;54(3):411--424.

\bibitem{Margolin2002}
Margolin LG, Shashkov M. {\it Second-order sign-preserving remapping on general
  grids. } Technical Report LA-UR-02-525: Los {A}lamos {N}ational {L}aboratory;
  2002.

\bibitem{Mavriplis2003}
Mavriplis DJ. Revisiting the least-squares procedure for gradient
  reconstruction on unstructured meshes.  In: AIAA 2003-3986. 16th {AIAA}
  {C}omputational {F}luid {D}ynamics {C}onference; June 23--26, 2003; Orlando,
  {F}lorida.

\bibitem{Scovazzi2008}
Scovazzi G, Love E, Shashkov M. Multi-scale {L}agrangian shock hydrodynamics on
  {Q1/P0} finite elements: {T}heoretical framework and two-dimensional
  computations.  {\it Comput {M}ethod {A}ppl {M}. }2008;197(9--12):1056--1079.

\end{thebibliography}


\begin{thebibliography}{10}
\providecommand \doibase [0]{http://dx.doi.org/}%

\bibitem{hughes1988space}
Hughes TJ, Hulbert GM. Space-time finite element methods for elastodynamics:
  Formulations and error estimates. {\it Computer Methods in Applied Mechanics
  and Engineering} 1988\string; 66(3)\string: 339-363.
\newblock \href {\doibase https://doi.org/10.1016/0045-7825(88)90006-0} {doi:
  https://doi.org/10.1016/0045-7825(88)90006-0}

\bibitem{langer2019coeff}
Langer U, Neum{\"u}ller M, Schafelner A. Space-Time Finite Element Methods for
  Parabolic Evolution Problems with Variable Coefficients. In:  Apel T, Langer
  U, Meyer A, Steinbach O. \kern-2pt, eds. {\it Advanced Finite Element Methods
  with Applications: Selected Papers from the 30th Chemnitz Finite Element
  Symposium 2017$\;$}Lecture Notes in Computational Science and Engineering.
  Springer; 2019\string: 247--275

\bibitem{sivas2021air}
Sivas AA, Southworth BS, Rhebergen S. AIR Algebraic Multigrid for a Space-Time
  Hybridizable Discontinuous Galerkin Discretization of Advection(-Diffusion).
  {\it SIAM Journal on Scientific Computing} 2021\string; 43(5)\string:
  A3393-A3416.
\newblock \href {\doibase 10.1137/20M1375103} {doi: 10.1137/20M1375103}

\bibitem{tezduyar1992new}
Tezduyar TE, Behr M, Liou J. A new strategy for finite element computations
  involving moving boundaries and interfaces--the DSD/ST procedure: I. The
  concept and the preliminary numerical tests. {\it Computer Methods in Applied
  Mechanics and Engineering} 1992\string; 94(3)\string: 339--351.

\bibitem{hubner2004monolithic}
H{\"u}bner B, Walhorn E, Dinkler D. A monolithic approach to fluid--structure
  interaction using space--time finite elements. {\it Computer Methods in
  Applied Mechanics and Engineering} 2004\string; 193(23-26)\string:
  2087--2104.

\bibitem{sathe2008modeling}
Sathe S, Tezduyar TE. Modeling of fluid--structure interactions with the
  space--time finite elements: contact problems. {\it Computational Mechanics}
  2008\string; 43(1)\string: 51.

\bibitem{spenke2020multi}
Spenke T, Hosters N, Behr M. A multi-vector interface quasi-Newton method with
  linear complexity for partitioned fluid--structure interaction. {\it Computer
  Methods in Applied Mechanics and Engineering} 2020\string; 361\string:
  112810.

\bibitem{behr2008simplex}
Behr M. Simplex space-time meshes in finite element simulations. {\it
  International Journal for Numerical Methods in Fluids} 2008\string;
  57(9)\string: 1421--1434.
\newblock \href {\doibase 10.1002/fld.1796} {doi: 10.1002/fld.1796}

\bibitem{danwitz2021four}
{von}~Danwitz M, Antony P, Key F, Hosters N, Behr M. Four-dimensional
  elastically deformed simplex space-time meshes for domains with time-variant
  topology. {\it International Journal for Numerical Methods in Fluids}
  2021\string; 93(12)\string: 3490-3506.
\newblock \href {\doibase https://doi.org/10.1002/fld.5042} {doi:
  https://doi.org/10.1002/fld.5042}

\bibitem{karyofylli2018simplex}
Karyofylli V, Frings M, Elgeti S, Behr M. Simplex space-time meshes in
  two-phase flow simulations. {\it International Journal for Numerical Methods
  in Fluids} 2018\string; 86\string: 218-230.
\newblock \href {\doibase 10.1002/fld.4414} {doi: 10.1002/fld.4414}

\bibitem{karyofylli2019simplex}
Karyofylli V, Wendling L, Make M, Hosters N, Behr M. Simplex space-time meshes
  in thermally coupled two-phase flow simulations of mold filling. {\it
  Computers \& Fluids} 2019\string; 192\string: 104261.

\bibitem{gesenhues2021simulating}
Gesenhues L, Behr M. Simulating dense granular flow using the $\mu$(I)-rheology
  within a space-time framework. {\it International Journal for Numerical
  Methods in Fluids} 2021\string; 93(9)\string: 2889-2904.
\newblock \href {\doibase https://doi.org/10.1002/fld.5014} {doi:
  https://doi.org/10.1002/fld.5014}

\bibitem{rendall2012conservative}
Rendall TC, Allen CB, Power ED. Conservative unsteady aerodynamic simulation of
  arbitrary boundary motion using structured and unstructured meshes in time.
  {\it International Journal for Numerical Methods in Fluids} 2012\string;
  70(12)\string: 1518--1542.

\bibitem{wang2015high}
Wang L, Persson PO. A high-order discontinuous Galerkin method with
  unstructured space--time meshes for two-dimensional compressible flows on
  domains with large deformations. {\it Computers \& Fluids} 2015\string;
  118\string: 53--68.

\bibitem{danwitz2019simplex}
{von Danwitz} M, Karyofylli V, Hosters N, Behr M. Simplex space-time meshes in
  compressible flow simulations. {\it International Journal for Numerical
  Methods in Fluids} 2019\string; 91(0)\string: 29-48.
\newblock \href {\doibase 10.1002/fld.4743} {doi: 10.1002/fld.4743}

\bibitem{karabelas2019generating}
Karabelas E, Neum{\"u}ller M. Generating admissible space-time meshes for
  moving domains in (d + 1) dimensions. In:  Langer U, Steinbach O. \kern-2pt,
  eds. {\it Space-Time Methods$\;$}Radon Series on Computational and Applied
  Mathematics. De Gruyter; 2019; Berlin\string: 185 -- 206.

\bibitem{boissonnat2021triangulating}
Boissonnat JD, Kachanovich S, Wintraecken M. Triangulating {Submanifolds}: {An}
  {Elementary} and {Quantified} {Version} of {Whitney}'s {Method}. {\it
  Discrete \& Computational Geometry} 2021\string; 66(1)\string: 386--434.
\newblock \href {\doibase 10.1007/s00454-020-00250-8} {doi:
  10.1007/s00454-020-00250-8}

\bibitem{caplan2020four}
Caplan PC, Haimes R, Darmofal DL, Galbraith MC. Four-Dimensional Anisotropic
  Mesh Adaptation. {\it Computer-Aided Design} 2020\string; 129\string: 102915.

\bibitem{takizawa2019node}
Takizawa K, Ueda Y, Tezduyar TE. A node-numbering-invariant directional length
  scale for simplex elements. {\it Mathematical Models and Methods in Applied
  Sciences} 2019\string; 29(14)\string: 2719-2753.

\bibitem{frontin2021foundations}
Frontin CV, Walters GS, Witherden FD, Lee CW, Williams DM, Darmofal DL.
  Foundations of space-time finite element methods: Polytopes, interpolation,
  and integration. {\it Applied Numerical Mathematics} 2021\string; 166\string:
  92-113.
\newblock \href {\doibase https://doi.org/10.1016/j.apnum.2021.03.019} {doi:
  https://doi.org/10.1016/j.apnum.2021.03.019}

\bibitem{elman2014finite}
Elman HC, Silvester DJ, Wathen AJ. {\it Finite elements and fast iterative
  solvers: with applications in incompressible fluid dynamics}.
\newblock Oxford University Press.
\newblock 2~ed. 2014.

\bibitem{shakib1991new}
Shakib F, Hughes TJ. A new finite element formulation for computational fluid
  dynamics: IX. Fourier analysis of space-time Galerkin/least-squares
  algorithms. {\it Computer Methods in Applied Mechanics and Engineering}
  1991\string; 87(1)\string: 35 - 58.
\newblock \href {\doibase 10.1016/0045-7825(91)90145-v} {doi:
  10.1016/0045-7825(91)90145-v}

\bibitem{donea2003finite}
Donea J, Huerta A. {\it Finite element methods for flow problems}.
\newblock John Wiley \& Sons.
\newblock 1~ed. 2003.

\bibitem{lozinski2009anisotropic}
Lozinski A, Picasso M, Prachittham V. An Anisotropic Error Estimator For The
  Crank-Nicolson Method: Application To A Parabolic Problem. {\it Siam Journal
  On Scientific Computing} 2009\string; 31\string: 2757-2783.
\newblock \href {\doibase 10.1137/080715135} {doi: 10.1137/080715135}

\bibitem{dubuis2018adaptive}
Dubuis S, Picasso M. An Adaptive Algorithm for the Time Dependent Transport
  Equation with Anisotropic Finite Elements and the Crank--Nicolson Scheme.
  {\it Journal of Scientific Computing} 2018\string; 75(1)\string: 350--375.
\newblock \href {\doibase 10.1007/s10915-017-0537-1} {doi:
  10.1007/s10915-017-0537-1}

\bibitem{hansbo2000crank}
Hansbo P. A Crank--Nicolson Type Space--Time Finite Element Method for
  Computing on Moving Meshes. {\it Journal of Computational Physics}
  2000\string; 159(2)\string: 274 - 289.
\newblock \href {\doibase 10.1006/jcph.2000.6436} {doi: 10.1006/jcph.2000.6436}

\bibitem{aziz1989continuous}
Aziz AK, Monk P. Continuous finite elements in space and time for the heat
  equation. {\it Mathematics of Computation} 1989\string; 52(186)\string:
  255--274.
\newblock \href {\doibase 10.1090/s0025-5718-1989-0983310-2} {doi:
  10.1090/s0025-5718-1989-0983310-2}

\bibitem{steinbach2015space}
Steinbach O. Space-time finite element methods for parabolic problems. {\it
  Computational methods in applied mathematics} 2015\string; 15(4)\string:
  551--566.
\newblock \href {\doibase 10.1515/cmam-2015-0026} {doi: 10.1515/cmam-2015-0026}

\bibitem{langer2019spacetime}
Langer U, Schafelner A. Space-Time Finite Element Methods for Parabolic
  Initial-Boundary Value Problems with Non-smooth Solutions. In: Springer. ;
  2020\string: 593--600

\bibitem{langer2021hex}
Langer U, Schafelner A. Space-time hexahedral finite element methods for
  parabolic evolution problems. arXiv;  2021.

\bibitem{langer2020efficient}
Langer U, Zank M. Efficient Direct Space-Time Finite Element Solvers for
  Parabolic Initial-Boundary Value Problems in Anisotropic Sobolev Spaces. {\it
  SIAM Journal on Scientific Computing} 2021\string; 43.
\newblock \href {\doibase 10.1137/20m1358128} {doi: 10.1137/20m1358128}

\bibitem{voulis2019discontinuous}
Voulis I, Reusken A. Discontinuous Galerkin time discretization methods for
  parabolic problems with linear constraints. {\it Journal of Numerical
  Mathematics} 2019\string; 27(3)\string: 155--182.
\newblock \href {\doibase 10.1515/jnma-2018-0013} {doi: 10.1515/jnma-2018-0013}

\bibitem{quarteroni2006numerical}
Quarteroni A, Sacco R, Saleri F. {\it Numerical Mathematics}.
\newblock Springer.
\newblock 2~ed. 2006.

\bibitem{ern2013theory}
Ern A, Guermond JL. {\it Theory and practice of finite elements}.
\newblock Springer.
\newblock 1~ed. 2004.

\bibitem{brooks1982streamline}
Brooks AN, Hughes TJR. Streamline upwind/Petrov-Galerkin formulations for
  convection dominated flows with particular emphasis on the incompressible
  Navier-Stokes equations. {\it Computer Methods in Applied Mechanics and
  Engineering} 1982\string; 32(1-3)\string: 199--259.

\bibitem{knechtges2018simulation}
Knechtges P. {\it Simulation of Viscoelastic Free-Surface Flows}. PhD thesis.
  RWTH Aachen University, Aachen;  2018.

\bibitem{jansen1999better}
Jansen KE, Collis SS, Whiting C, Shakib F. A better consistency for low-order
  stabilized finite element methods. {\it Computer Methods in Applied Mechanics
  and Engineering} 1999\string; 174(1)\string: 153 - 170.
\newblock \href {\doibase 10.1016/s0045-7825(98)00284-9} {doi:
  10.1016/s0045-7825(98)00284-9}

\bibitem{mojtabi2015one}
Mojtabi A, Deville MO. One-dimensional linear advection--diffusion equation:
  Analytical and finite element solutions. {\it Computers \& Fluids}
  2015\string; 107\string: 189--195.
\newblock \href {\doibase 10.1016/j.compfluid.2014.11.006} {doi:
  10.1016/j.compfluid.2014.11.006}

\bibitem{thomee2006Galerkin}
Thom{\'e}e V. {\it Galerkin Finite Element Methods for Parabolic Problems}.
\newblock Springer.
\newblock 2~ed. 2006.

\bibitem{demkowicz1986adaptive}
Demkowicz L, Oden J. An adaptive characteristic Petrov-Galerkin finite element
  method for convection-dominated linear and nonlinear parabolic problems in
  one space variable. {\it Journal of Computational Physics} 1986\string;
  67(1)\string: 188--213.
\newblock \href {\doibase 10.1016/0021-9991(86)90121-x} {doi:
  10.1016/0021-9991(86)90121-x}

\bibitem{geuzaine2009gmsh}
Geuzaine C, Remacle JF. Gmsh: A 3-D finite element mesh generator with built-in
  pre-and post-processing facilities. {\it International Journal for Numerical
  Methods in Engineering} 2009\string; 79(11)\string: 1309--1331.

\end{thebibliography}

\newpage
\appendix
\section{Appendix}
\label{sec:appendix}
\renewcommand{\plthgtTD}{0.33\textheight}
\renewcommand{\pltwdtTD}{0.5\textwidth}

\newcommand{\logPlotThreeDApp}[1]{
{\footnotesize
%\tikzsetnextfilename{app-#1}
\begin{tikzpicture}
\begin{axis}
[
height=\plthgtTD, width=\pltwdtTD,
point meta min = -9,
zmax=0, 
zmin=-12,
unbounded coords=jump,
y dir=reverse,
xlabel=$l$, ylabel=$m$, 
zlabel=$\log (E_{lm})$,
z label style={at={(axis description cs:0.0,0.95)},,rotate=-90,anchor=north},
]
\addplot3 [surf, mesh/rows=15, mesh/ordering=y varies, z filter/.expression= {z==-15 ? nan : z}] table [col sep=comma]  {Data/#1L2ErrorCoords.csv};
\end{axis}
\end{tikzpicture}
}
}

\newcommand{\sfl}{.48\textwidth}

\begin{figure}[h]
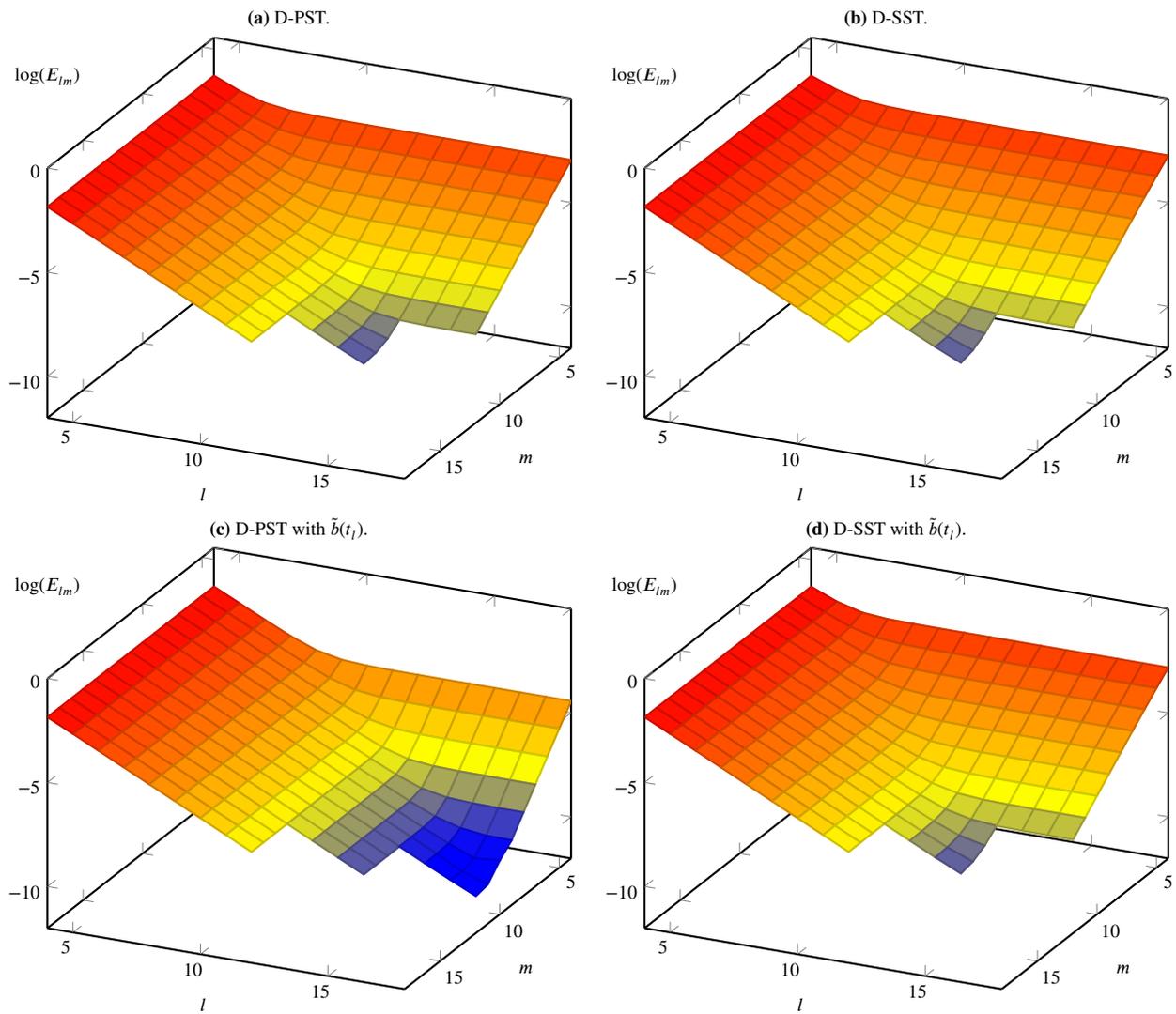

\centering
%	\begin{subfigure}{\sfl}
%		\centering
%		\caption{D-PST.}
%		\label{fig:appFSTP2}
%		\vspace{1ex}
%		\logPlotThreeDApp{FSTQuadNu0.1P2}
%	\end{subfigure}
%	\begin{subfigure}{\sfl}
%		\centering
%		\caption{D-SST.}
%		\label{fig:appSSTP2}
%		\vspace{1ex}
%		\logPlotThreeDApp{SST1DNu0.1P2}
%	\end{subfigure}
%	\begin{subfigure}{\sfl}
%		\centering
%		\caption{D-PST with $\tilde{b}(t_l)$.}
%		\label{fig:appFSTP2Mod}
%		\vspace{1ex}
%		\logPlotThreeDApp{FSTQuadNu0.1P2Mod}
%	\end{subfigure}
%	\begin{subfigure}{\sfl}
%		\centering
%		\caption{D-SST with $\tilde{b}(t_l)$.}
%		\label{fig:appSSTP2Mod}
%		\vspace{1ex}
%		\logPlotThreeDApp{SST1DNu0.1P2Mod}
%	\end{subfigure}
\subcaptionbox{D-PST. \label{fig:appFSTP2}}{\logPlotThreeDApp{FSTQuadNu0.1P2}}%
\subcaptionbox{D-SST. \label{fig:appSSTP2}}{\logPlotThreeDApp{SST1DNu0.1P2}} \\ \vspace{1ex}
\subcaptionbox{D-PST with $\tilde{b}(t_l)$. \label{fig:appFSTP2Mod}}{\logPlotThreeDApp{FSTQuadNu0.1P2Mod}}%
\subcaptionbox{D-SST with $\tilde{b}(t_l)$. \label{fig:appSSTP2Mod}}{\logPlotThreeDApp{SST1DNu0.1P2Mod}} \\
\caption{Convergence of nodal error measure for parabolic problem IBVP 2.}
\label{fig:appC2Ne}
\end{figure} 

\newcommand{\logPlotThreeDAppReal}[1]{
{\footnotesize
%\tikzsetnextfilename{app-Real#1}
\begin{tikzpicture}
\begin{axis}
[
height=\plthgtTD, width=\pltwdtTD,
point meta min = -8,
zmax=0,
zmin=-12,
unbounded coords=jump,
y dir=reverse,
xlabel=$l$, ylabel=$m$, 
zlabel=$\log (e_{lm})$,
z label style={at={(axis description cs:0.0,0.95)},,rotate=-90,anchor=north},
]
\addplot3 [surf, mesh/rows=15, mesh/ordering=y varies, z filter/.expression= {z==-15 ? nan : z}] table [col sep=comma]  {Data/#1RealL2ErrorCoords.csv};
\end{axis}
\end{tikzpicture}
}
}

\begin{figure}
\centering
\subcaptionbox{D-PST, $\Peclet = 0$. \label{fig:appFSTPeZero}}{\logPlotThreeDAppReal{FSTQuadPe0}}%
\subcaptionbox{D-SST, $\Peclet = 0$. \label{fig:appSSTPeZero}}{\logPlotThreeDAppReal{SST1DPe0}} \\ \vspace{1ex}
\subcaptionbox{D-PST, $\Peclet = 1$. \label{fig:appFSTNuTen}}{\logPlotThreeDAppReal{FSTQuadNu0.1}}%
\subcaptionbox{D-SST, $\Peclet = 1$. \label{fig:appSSTNuTen}}{\logPlotThreeDAppReal{SST1DNu0.1}} \\ \vspace{1ex}
\subcaptionbox{D-PST, $\Peclet = 10$. \label{fig:appFSTNuHundred}}{\logPlotThreeDAppReal{FSTQuadNu0.01}}%
\subcaptionbox{D-SST, $\Peclet = 10$. \label{fig:appSSTNuHundred}}{\logPlotThreeDAppReal{SST1DNu0.01}} \\
\end{figure} 

\begin{figure}
\centering
\ContinuedFloat 
\subcaptionbox{D-PST, $\Peclet = 100$. \label{fig:appFSTNuThousand}}{\logPlotThreeDAppReal{FSTQuadNu0.001}}%
\subcaptionbox{D-SST, $\Peclet = 100$. \label{fig:appSSTNuThousand}}{\logPlotThreeDAppReal{SST1DNu0.001}} \\ \vspace{1ex}
\subcaptionbox{D-PST, $\Peclet = 1000$. \label{fig:appFSTNuTenThousand}}{\logPlotThreeDAppReal{FSTQuadNu0.0001}}%
\subcaptionbox{D-SST, $\Peclet = 1000$. \label{fig:appSSTNuTenThousand}}{\logPlotThreeDAppReal{SST1DNu0.0001}} \\ \vspace{1ex}
\subcaptionbox{D-PST,  $\Peclet = \infty$. \label{fig:appFSTNuZero}}{\logPlotThreeDAppReal{FSTQuadNu0.0}}%
\subcaptionbox{D-SST,  $\Peclet = \infty$. \label{fig:appSSTNuZero}}{\logPlotThreeDAppReal{SST1DNu0.0}} \\
\caption{Convergence visualization of $L_2$ error for model problem IBVP 1.}
\label{fig:appADC1L2}
\end{figure}

\begin{figure}
\centering
\subcaptionbox{D-PST, $\Peclet = 0$. \label{fig:appNeFSTPeZero}}{\logPlotThreeDApp{FSTQuadPe0}}%
\subcaptionbox{D-SST, $\Peclet = 0$. \label{fig:appNeSSTPeZero}}{\logPlotThreeDApp{SST1DPe0}} \\ \vspace{1ex}
\subcaptionbox{D-PST, $\Peclet = 1$. \label{fig:appNeFSTNuTen}}{\logPlotThreeDApp{FSTQuadNu0.1}}%
\subcaptionbox{D-SST, $\Peclet = 1$. \label{fig:appNeSSTNuTen}}{\logPlotThreeDApp{SST1DNu0.1}} \\ \vspace{1ex}
\subcaptionbox{D-PST, $\Peclet = 10$. \label{fig:appNeFSTNuHundred}}{\logPlotThreeDApp{FSTQuadNu0.01}}%
\subcaptionbox{D-SST, $\Peclet = 10$. \label{fig:appNeSSTNuHundred}}{\logPlotThreeDApp{SST1DNu0.01}} \\ 
\end{figure} 

\begin{figure}
\centering
\ContinuedFloat 
\subcaptionbox{D-PST, $\Peclet = 100$. \label{fig:appNeFSTNuThousand}}{\logPlotThreeDApp{FSTQuadNu0.001}}%
\subcaptionbox{D-SST, $\Peclet = 100$. \label{fig:appNeSSTNuThousand}}{\logPlotThreeDApp{SST1DNu0.001}} \\ \vspace{1ex}
\subcaptionbox{D-PST, $\Peclet = 1000$. \label{fig:appNeFSTNuTenThousand}}{\logPlotThreeDApp{FSTQuadNu0.0001}}%
\subcaptionbox{D-SST, $\Peclet = 1000$. \label{fig:appNeSSTNuTenThousand}}{\logPlotThreeDApp{SST1DNu0.0001}} \\ \vspace{1ex}
\subcaptionbox{D-PST, $\Peclet = \infty$. \label{fig:appNeFSTNuZero}}{\logPlotThreeDApp{FSTQuadNu0.0}}%
\subcaptionbox{D-SST, $\Peclet = \infty$. \label{fig:appNeSSTNuZero}}{\logPlotThreeDApp{SST1DNu0.0}} \\
\caption{Convergence visualization of nodal error for model problem IBVP 1.}
\label{fig:appADC1Ne}
\end{figure}

\end{document}